\documentclass[10pt,journal]{IEEEtran}
\usepackage{graphicx}	
\usepackage[cmex10]{amsmath}
\usepackage{amssymb}
\usepackage{amsthm}
\usepackage{bm}
\usepackage{color}
\usepackage{filecontents,lipsum}
\usepackage{paralist}
\usepackage{enumitem}
\usepackage{cite}
\usepackage{longtable}

\makeatletter
\g@addto@macro\normalsize{
  \setlength\abovedisplayskip{3pt}
  \setlength\belowdisplayskip{3pt}
  \setlength\abovedisplayshortskip{3pt}
  \setlength\belowdisplayshortskip{3pt}
}
\makeatother

\makeatletter
\def\thm@space@setup{
  \thm@preskip=3pt plus 1pt minus 1pt
  \thm@postskip=3pt plus 1pt minus 1pt
}
\makeatother

\theoremstyle{plain}
\newtheorem{thm}{Theorem}

\newtheorem{prop}{Proposition}

\newtheorem{assu}{Assumption}

\theoremstyle{definition}
\newtheorem{defn}{Definition}

\theoremstyle{remark}
\newtheorem{rem}{Remark}

\newcommand{\lc}{\left\{}

\newcommand{\rc}{\right\}}

\usepackage{amsfonts}

\title{{\huge A Polynomial-Time Method for Testing Admissibility of Uncertain Power Injections in Microgrids}}

\author{{\normalsize Cong~Wang,~\IEEEmembership{Student Member,~IEEE}, Eleni~Stai,~\IEEEmembership{Member,~IEEE}, and Jean-Yves~Le~Boudec,~\IEEEmembership{Fellow,~IEEE}}

\thanks{This work is supported by SNSF-NRP70 ``Energy Turnaround''.

The authors are with \'{E}cole Polytechnique F\'{e}d\'{e}rale de Lausanne (EPFL), CH-1015 Lausanne, Switzerland (e-mail: cong.wang@epfl.ch; eleni.stai@epfl.ch; jean-yves.leboudec@epfl.ch).}}

\begin{document}
\maketitle

\begin{abstract}
We study the admissibility of power injections in single-phase microgrids, where the electrical state is represented by complex nodal voltages and controlled by nodal power injections. Assume that (i) there is an initial electrical state that satisfies security constraints and the non-singularity of load-flow Jacobian, and (ii) power injections reside in some uncertainty set. We say that the uncertainty set is admissible for the initial electrical state if any continuous trajectory of the electrical state is ensured to be secured and non-singular as long as power injections remain in the uncertainty set. We use the recently proposed V-control and show two new results.
 First, if a complex nodal voltage set V is convex and every element in V is non-singular, then V is a domain of uniqueness. Second, we give sufficient conditions to guarantee that every element in some power injection set S has a load-flow solution in V, based on impossibility of obtaining load-flow solutions at the boundary of V. By these results, we develop a framework for the admissibility-test method; this framework is extensible to multi-phase grids. Within the framework, we establish a polynomial-time method, using the infeasibility check of convex optimizations. The method is evaluated numerically.
\end{abstract}

{
\section*{Nomenclature}
\begin{small}
\hspace{-0.4cm}
\begin{tabular}{ll}
$N\qquad\qquad\qquad\qquad\quad$ & Number of $PQ$ buses \\
$\mathcal{N}=\{0,...,N\}$ & Set of buses, 0 for the slack bus \\
$\mathcal{N}^{PQ}=\mathcal{N}\setminus\{0\}$ \\
$\mathcal{E}$ & Set of ordered index pairs for \\
& referring to specific branch current \\
$v_j$ & Complex nodal voltage at bus $j\in\mathcal{N}$ \\
$\mathbf{v}=(v_1,...,v_N)^T$ & Complex nodal voltage vector \\
$\mathbf{v}^\mathrm{initial}$ & Initial complex nodal voltage vector \\
$\mathbf{w}$ & Zero-load complex nodal voltage vector \\
$i_{jk}$ & Complex branch current from bus \\
& $j$ to $k$, for ${jk}\in\mathcal{E}$ \\
$i_j$ & Complex nodal current at bus $j\in\mathcal{N}$ \\
$\mathbf{i}=(i_1,...,i_N)^T$ & Complex nodal current vector \\
$s_j$ & Complex nodal power injection \\
& at bus $j\in\mathcal{N}$ \\
$\mathbf{s}=(s_1,...,s_N)^T$ & Complex nodal power injection vector \\
$\mathbf{F}()$ & Function on $\mathbb{C}^{N}$ that maps any $\mathbf{v}$ \\
& into its corresponding $\mathbf{s}$ \\
$\mathbf{J}_\mathbf{F}(\mathbf{v})$ & Jacobian of $\mathbf{F}$ at $\mathbf{v}$ \\
$\mathcal{V}$ & Set of $\mathbf{v}$ \\
$\mathcal{S}$ & Set of $\mathbf{s}$ \\
$\mathcal{S}^\mathrm{uncertain}$ & Uncertainty set of $\mathbf{s}$ \\
$\mathbf{Y}$ & Nodal admittance matrix \\
$\mathbf{Y}_{LL}$ & Submatrix of $\mathbf{Y}$ \\
$V^\mathrm{min}_j$, $V^\mathrm{max}_j$, and $I^\mathrm{max}_{jk}$ & Security bounds for $v_j,j\in\mathcal{N}^{PQ}$ \\
& and $i_{jk},jk\in\mathcal{E}$ \\
$\mathrm{Re}()$ and $\mathrm{Im}()$ & Real and imaginary parts of a \\
& complex variable \\
\end{tabular}
\end{small}
\begin{small}
\begin{tabular}{ll}
$f_j^\mathrm{V,low}$, $f_j^\mathrm{V,up}$, $f_{jk}^\mathrm{I,branch}$ & Polynomials of $\mathrm{Re}(\mathbf{v})$, $\mathrm{Im}(\mathbf{v})$ that are \\
& used to express security constraints, \\
& where $j\in\mathcal{N}^{PQ}$ and $jk\in\mathcal{E}$ \\
$\tilde{f}_{jk}^\mathrm{I,branch}$ and $\tilde{f}_j^\mathrm{I,node}$ & Polynomials of $\mathrm{Re}(\mathbf{v})$, $\mathrm{Im}(\mathbf{v})$ that are \\
& used in the proposed method, \\
& where $j\in\mathcal{N}^{PQ}$ and $jk\in\mathcal{E}$ \\
$~\bar{\cdot}~$ & Complex conjugation \\
$\|\cdot\|_1$ and $\|\cdot\|_\infty$ & $\ell_1$ and $\ell_\infty$ norm \\
$(\cdot)_{m,n}$ and $(\cdot)_n$ & Respectively the entry of $m$-th row, \\
& $n$-th column in a matrix, and \\
& the $n$-th entry in a vector \\
$\mathrm{Row}_j()$ & The $j$-th row of a matrix \\
$\partial\mathcal{V}$ & Topological boundary of set $\mathcal{V}$ \\
$\kappa$ & Real scaling factor
\end{tabular}
\end{small}
}

\section{Introduction}
\subsection{Background}
In the last decade, there has been a large number of excellent works on microgrids. These works range from theoretical aspects to real-world applications (see \cite{MicroCostMin,Micro2,Micro3,Micro4,Micro5,Micro6,Micro7} for some examples). A prominent feature of many microgrids is the integration of renewable energy sources and electrical vehicles. Compared to the fixed power generation and consumption in traditional power grids, the nodal power injections in modern microgrids are usually uncertain due to the volatility of these sources and loads. In practice, such uncertain nodal power injections might result in some undesired electrical state, when the latter is controlled by the former.

For example, consider a microgrid control system that computes setpoints of nodal power injections and sends them as explicit commands to the grid resources for implementation. Due to aforementioned volatility, the control system cannot be sure that these setpoints of nodal power injections will be exactly implemented but knows that they will reside in some uncertainty set {\small $\mathcal{S}^\mathrm{uncertain}$} \cite{Robust,commelec1}. Now, assume that (i) the electrical state is represented by the steady-state complex nodal voltages, (ii) the initial electrical state fulfills a specific set of security constraints and long-term voltage stability (i.e., the load-flow Jacobian is non-singular), and that (iii) the initial nodal power injections are included in {\small $\mathcal{S}^\mathrm{uncertain}$}. For the control system, it would like to be sure that the electrical state continues to satisfy the security constraints and long-term voltage stability, as long as the implemented nodal power injections stay in {\small $\mathcal{S}^\mathrm{uncertain}$}. Obviously, for this to occur, every element in {\small $\mathcal{S}^\mathrm{uncertain}$} needs to have at least one load-flow solution that is secured and non-singular. However, as pointed out in \cite{Vcontrol}, this is not sufficient.

The above example explains the formulation of the ``admissibility problem'' in this paper. More precisely, given a secured and non-singular initial electrical state as well as an uncertainty set {\small $\mathcal{S}^\mathrm{uncertain}$} that contains the initial nodal power injections, we say that {\small $\mathcal{S}^\mathrm{uncertain}$} is admissible for this initial electrical state if any continuous trajectory of the electrical state is ensured to be secured and non-singular as long as the corresponding nodal power injections are in {\small $\mathcal{S}^\mathrm{uncertain}$} (a formal definition is presented in Section \ref{sec:pbf}). Here, by ``continuous trajectory'', we mean that the electrical state (represented by the steady-state complex nodal voltages) changes as a continuous function of time.

To solve the admissibility problem, we propose to use the recently developed theory of {\small $\mathcal{V}$}-control \cite{Vcontrol}. In short, let the electrical state be represented by complex nodal voltages and {\small $\mathcal{V}$} be a set of complex nodal voltages. Then, an arbitrary set {\small $\mathcal{S}$} of complex nodal power injections is a ``domain of {\small $\mathcal{V}$}-control'' if any continuous trajectory of the electrical state that starts in {\small $\mathcal{V}$} must stay in {\small $\mathcal{V}$}, as long as the corresponding trajectory of nodal power injections stays in {\small $\mathcal{S}$}.

A naive application of the theory of {\small $\mathcal{V}$}-control to the admissibility problem would be showing that {\small $\mathcal{S}^\mathrm{uncertain}$} is a domain of {\small $\mathcal{V}$}-control, where {\small $\mathcal{V}$} is the set of all secured and non-singular electrical states. By Lemma 2 in \cite{Vcontrol} (recalled in this paper as Theorem~\ref{thm:AA}), under the assumption that all security constraints are strict inequalities, a sufficient condition for this {\small $\mathcal{V}$}-control to hold would be that every element in {\small $\mathcal{S}^\mathrm{uncertain}$} has exactly one corresponding electrical state in {\small $\mathcal{V}$}. However, this could be impractical, as the condition that every element in {\small $\mathcal{S}^\mathrm{uncertain}$} corresponds to a unique secured and non-singular electrical state might not hold. We demonstrate this in the next example.

\begin{figure}[t]
\begin{center}
\includegraphics[scale=0.66]{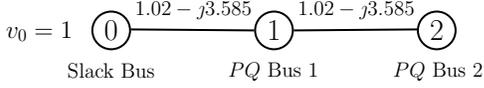}
\caption{Grid topology with slack-bus voltage and series admittances in p.u.}
\label{fig:ThreeBusNetwork}
\end{center}
\end{figure}

Consider the grid in Figure \ref{fig:ThreeBusNetwork}, which is formed by three serially connected buses \cite{UnControl1}. For ease of exposition, let the complex nodal voltage and the complex nodal power injection at {\small $PQ$} bus {\small $j\in\{1,2\}$} be {\small $v_j$}, {\small $s_j$}, respectively. Now, suppose that

\begin{itemize}
\item An electrical state is secured if the deviations in nodal voltage magnitudes are less than {\small $\pm10\%$} of the nominal value (here, security bounds on branch current magnitudes are chosen to be sufficiently large hence do not come into effect);
\item The initial nodal power injections are {\small $s_1^\mathrm{initial}=-1.105+\jmath1$}, {\small $s_2^\mathrm{initial}=-1+\jmath1.105$} in p.u.;
\item {\small $\mathcal{S}^\mathrm{uncertain}=\{(s_1;s_2):|s_j-s_j^\mathrm{initial}|\leq10^{-5},j\in\{1,2\}\}.$}
\end{itemize}

\begin{figure}[t]
\begin{center}
\includegraphics[scale=0.47]{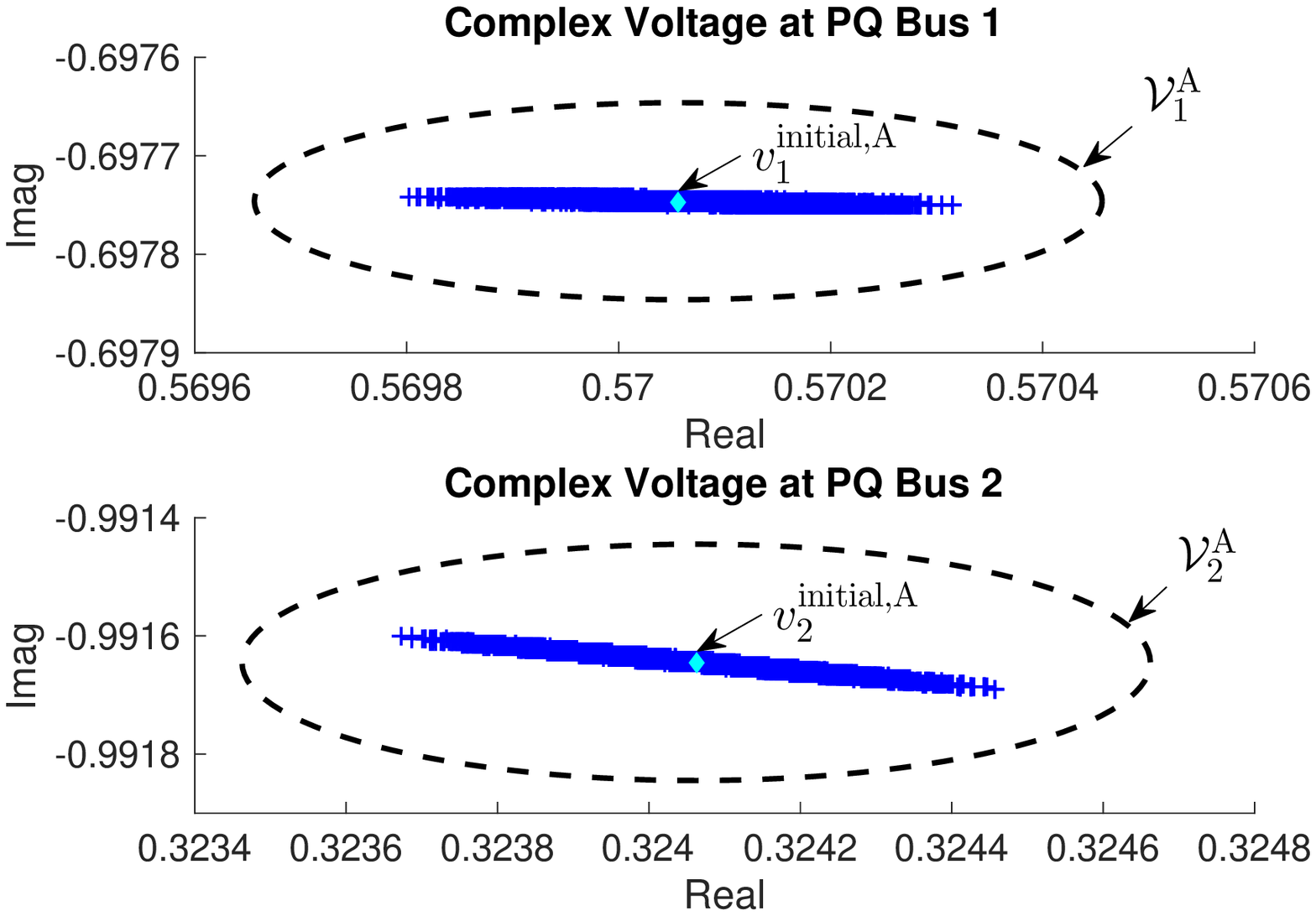}

(a)\\
\includegraphics[scale=0.47]{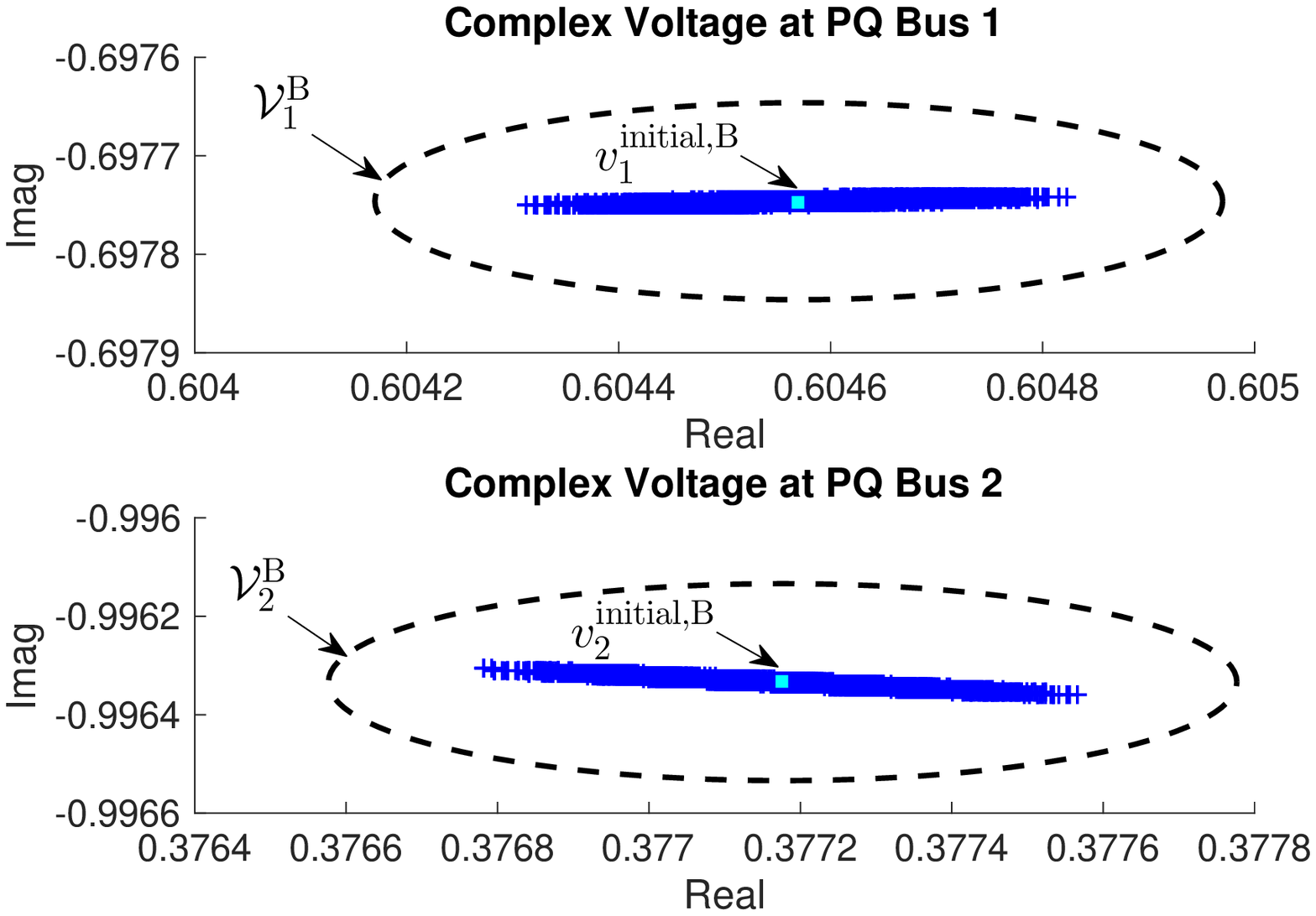}

(b)
\caption{(i) Candidate initial electrical states $(v_1^\mathrm{initial,A};v_2^\mathrm{initial,A})$ and $(v_1^\mathrm{initial,B};v_2^\mathrm{initial,B})$, which are presented by cyan ``diamond'' and ``square''; (ii) All secured and non-singular electrical states (marked by scattered blue ``$+$'') that correspond to the elements in $\mathcal{S}^\mathrm{uncertain}$; (iii) Open, secured and non-singular sets $\mathcal{V}_1^\mathrm{A}\times\mathcal{V}_2^\mathrm{A}$ and $\mathcal{V}_1^\mathrm{B}\times\mathcal{V}_2^\mathrm{B}$.}
\label{fig:ExampleCase}
\end{center}
\end{figure}

For the given initial nodal power injections, there are two corresponding secured and non-singular electrical states (as shown in Figure \ref{fig:ExampleCase}). These two candidate initial electrical states are denoted by {\small $(v_1^\mathrm{initial,A};v_2^\mathrm{initial,A})$} and {\small $(v_1^\mathrm{initial,B};v_2^\mathrm{initial,B})$}. Additionally, each element in {\small $\mathcal{S}^\mathrm{uncertain}$} has two corresponding secured and non-singular electrical states that are located around the two candidate initial electrical states. We numerically find these electrical states for all elements in {\small $\mathcal{S}^\mathrm{uncertain}$}, and we mark them by scattered ``$+$'' in Figure \ref{fig:ExampleCase}.

In this grid, {\small $\mathcal{S}^\mathrm{uncertain}$} is admissible for {\small $(v_1^\mathrm{initial,A};v_2^\mathrm{initial,A})$}. To see why, construct the open, secured and non-singular set (i.e., every electrical state in this set is secured and non-singular) {\small $\mathcal{V}^\mathrm{A}=\mathcal{V}_1^\mathrm{A}\times\mathcal{V}_2^\mathrm{A}$} as in Figure ~\ref{fig:ExampleCase}(a) with
\begin{itemize}
\item {\small $\mathcal{V}_1^\mathrm{A}=\{v_1:(\frac{\mathrm{Re}(v_1-v_1^\mathrm{initial,A})}{0.0004})^2+(\frac{\mathrm{Im}(v_1-v_1^\mathrm{initial,A})}{0.0001})^2<1\}$};
\item {\small $\mathcal{V}_2^\mathrm{A}=\{v_2:(\frac{\mathrm{Re}(v_2-v_2^\mathrm{initial,A})}{0.0006})^2+(\frac{\mathrm{Im}(v_2-v_2^\mathrm{initial,A})}{0.0002})^2<1\}$}.
\end{itemize}
Here, {\small $\mathrm{Re}()$} (resp. {\small $\mathrm{Im}()$}) is the real (resp. imaginary) part of a complex number, and ``{\small $\times$}'' means the Cartesian product. Clearly, for each element in {\small $\mathcal{S}^\mathrm{uncertain}$}, the corresponding electrical state around {\small $(v_1^\mathrm{initial,A};v_2^\mathrm{initial,A})$} is in {\small $\mathcal{V}^\mathrm{A}$}. Moreover, for each element in {\small $\mathcal{S}^\mathrm{uncertain}$}, there is exactly one corresponding electrical state in {\small $\mathcal{V}^\mathrm{A}$}. Next, by Lemma 2 in \cite{Vcontrol}, {\small $\mathcal{S}^\mathrm{uncertain}$} is a domain of {\small $\mathcal{V}^\mathrm{A}$}-control and is thus admissible for {\small $(v_1^\mathrm{initial,A};v_2^\mathrm{initial,A})$}. Note that the same logic can be employed to prove that {\small $\mathcal{S}^\mathrm{uncertain}$} is also admissible for {\small $(v_1^\mathrm{initial,B};v_2^\mathrm{initial,B})$}, using {\small $\mathcal{V}^\mathrm{B}=\mathcal{V}_1^\mathrm{B}\times\mathcal{V}_2^\mathrm{B}$} shown in Figure \ref{fig:ExampleCase}(b).

This example illustrates that, in order to apply the theory of {\small $\mathcal{V}$}-control in \cite{Vcontrol} to the admissibility problem, it is necessary to find an appropriate domain {\small $\mathcal{V}$} that is typically smaller than the set of all secured and non-singular electrical states. Using these observations of {\small $\mathcal{V}$}-control, we proposed in \cite{Vcontrol} a heuristic method for solving the admissibility problem. Although this heuristic method works in both single-phase and multi-phase grids, its performance can be unstable due to the exponential-time complexity. In order to address this issue, we develop a polynomial-time method in this paper.

\subsection{Main Contributions}
Our main contributions are as follows.

1) In Section \ref{sec:frame}, we propose a framework that, based on the theory of {\small $\mathcal{V}$}-control in \cite{Vcontrol}, solves the admissibility problem. The formal definition of {\small $\mathcal{V}$}-control is recalled in Definition \ref{defn:I} and some sufficient conditions for {\small $\mathcal{V}$}-control are recalled in Theorem \ref{thm:AA}. To satisfy these conditions, we propose two new theorems:
\begin{itemize}
\item In Theorem \ref{thm:BB}, we show that if a set {\small $\mathcal{V}$} of complex nodal voltages is convex and all its elements satisfy the non-singularity of the load-flow Jacobian, then it is a domain of uniqueness (i.e., if two elements in {\small $\mathcal{V}$} have the same image in power space under the load-flow function, then they must be identical).
\item In Theorem \ref{thm:CC}, we provide sufficient conditions on an arbitrary nodal voltage set {\small $\mathcal{V}$} and an arbitrary nodal power injection set {\small $\mathcal{S}$} so that every element in {\small $\mathcal{S}$} is guaranteed to have a load-flow solution in {\small $\mathcal{V}$}; the conditions are based on the impossibility of obtaining load-flow solutions at the boundary of {\small $\mathcal{V}$}.
\end{itemize}
Although the framework and the theorems are proposed for single-phase grids, they can be easily extended to multi-phase grids. With the framework, we can develop methods that are alternatives to the heuristic method in \cite{Vcontrol}.

 2) In Section \ref{sec:method}, we develop a concrete method that fits into the proposed framework and implements the aforementioned theorems. Compared to the heuristic method in \cite{Vcontrol}, a prominent feature of the method in this paper is the polynomial-time complexity. This feature is achieved by checking the infeasibility of convex optimizations, which is not considered in the heuristic method in \cite{Vcontrol}. To evaluate the performance of our polynomial-time method, we apply it to a few test grids. For illustration purposes, we show numerical results on (i) a meshed grid, (ii) the modified IEEE 13-Bus Test Feeder, and on (iii) the CIGRE North American LV Distribution Network in Section \ref{sec:numev}.

\section{Grid Model and Theoretical Foundations} \label{sec:modeldef}
\subsection{Grid Model}
The following developments rely on a quasi-stationary representation of the AC power grids; it models the grid parameters (impedances and admittances), state variables (complex nodal voltages, nodal currents, and branch currents), and forcing variables (complex nodal power injections) by means of phasors. In particular, inrushes are approximated by continuous power injections that are characterized with adequate ramp rates. Moreover, we assume that

\begin{itemize}
\item The method presented in the paper is integrated in control systems that are responsible for steering a grid. An example of such control systems is the Commelec \cite{commelec1};
\item During the application of the proposed method, the grid topology is known and does not change. This is, in practice, true in many systems. For instance, the topology of a European distribution grid changes only several times in a month due to line faults or specific operational plans.
\end{itemize}

Now, consider a single-phase grid with one slack bus, {\small $N$ $PQ$} buses and a generic topology (i.e., radial or meshed). \footnote{In \cite{Converter1,Converter3,Converter4,Converter5,Converter8}, it has been thoroughly discussed that buses with energy sources such as photovoltaic panels, wind turbines, microturbines, and fuel cells can be operated in $PQ$ mode via power electronic inverters.} In the paper, we assign index {\small $0$} to the slack bus, and indexes {\small $1,...,N$} to the {\small $PQ$} buses. For convenience of expression, we define
\begin{itemize}
\item {\small $\mathcal{N}\triangleq \{0,...,N\}$} as the index set of all buses;
\item {\small $\mathcal{N}^{PQ}\triangleq\mathcal{N}\setminus\{0\}$} as the index set of {\small $PQ$} buses;
\item {\small $\mathcal{E}\triangleq\{jk:$} a branch exists between buses {\small$j,k\in\mathcal{N}\}$} as a set of ordered index pairs for referring to branch currents.
\end{itemize}
At each bus {\small $j\in\mathcal{N}$}, we denote the complex phase-to-ground nodal voltage, nodal current and nodal power injection by {\small $v_j$}, {\small $i_j$} and {\small $s_j$}, respectively. Furthermore, let {\small $\mathbf{v}\triangleq(v_1,...,v_N)^T\in\mathbb{C}^N$}, {\small $\mathbf{i}\triangleq(i_1,...,i_N)^T\in\mathbb{C}^N$} and {\small $\mathbf{s}\triangleq(s_1,...,s_N)^T\in\mathbb{C}^N$}. We have that
\begin{itemize}
\item The branch current from bus {\small $j$} to {\small $k$} can be represented as
\begin{small}
\begin{equation} \label{eqn:bc}
i_{jk}=a_{jk}v_0+\mathbf{c}_{jk}^T\mathbf{v},
\end{equation}
\end{small}
where {\small $a_{jk}\in\mathbb{C},~\mathbf{c}_{jk}\in\mathbb{C}^N$} are constant and given by the passive transmission devices in \cite{PowerModel}. These passive devices include transmission lines (with shunt capacitance) and transformers. For example, let us think about a {\small $\pi$}-modeled transmission line between buses {\small $j$} and {\small $k$}. Suppose that the series admittance is {\small $y_{jk}^\mathrm{series}$} and the shunt element is {\small $b_{jk}^\mathrm{shunt}$}. Then, {\small $i_{jk}=y_{jk}^\mathrm{series}(v_j-v_k)+(\jmath b_{jk}^\mathrm{shunt}/2)v_j$}, i.e., equation \eqref{eqn:bc} holds.
\item Denote the {\small $(N+1)\times(N+1)$} nodal admittance matrix \cite{ComAdm} by {\small $\mathbf{Y}$}, which can be obtained using the grid topology and the passive transmission devices in \cite{PowerModel}. Then, {\small $\mathbf{v}$}, {\small $\mathbf{i}$}, {\small $\mathbf{s}$} fulfill the following load-flow equation system, where ``$~\bar{~}~$'' stands for complex conjugation.
\begin{small}
\begin{equation*}
\begin{bmatrix}
i_0 \\
\mathbf{i}
\end{bmatrix}=\mathbf{Y}
\begin{bmatrix}
v_0 \\
\mathbf{v}
\end{bmatrix},
\end{equation*}
\end{small}
\begin{small}
\begin{equation*}
\begin{bmatrix}
s_0 \\
\mathbf{s}
\end{bmatrix}=
\begin{bmatrix}
v_0 & \\
& \mathrm{diag}(\mathbf{v})
\end{bmatrix}
\begin{bmatrix}
\overline{i}_0 \\
\overline{\mathbf{i}}
\end{bmatrix}.
\end{equation*}
\end{small}
Note that {\small $\mathbf{Y}$} can be partitioned as 
\begin{small}
\begin{equation*}
\mathbf{Y}=
\begin{bmatrix}
\mathbf{Y}_{00} & \mathbf{Y}_{0L} \\
\mathbf{Y}_{L0} & \mathbf{Y}_{LL}
\end{bmatrix},
\end{equation*}
\end{small}
with {\small $N\times N$} matrix {\small $\mathbf{Y}_{LL}$} being invertible in practice \cite{EU}. In this way, by defining the vector of zero-load complex nodal voltages {\small $\mathbf{w}=-\mathbf{Y}^{-1}_{LL}\mathbf{Y}_{L0}v_0$}, we get that
\begin{small}
\begin{equation} \label{eqn:I}
\mathbf{i}=\mathbf{Y}_{LL}(\mathbf{v}-\mathbf{w}),
\end{equation}
\begin{align} \label{eqn:F}
\mathbf{s}&=\mathrm{diag}(\mathbf{v})\overline{\mathbf{Y}}_{LL}(\overline{\mathbf{v}}-\overline{\mathbf{w}}) \nonumber \\
&\triangleq\mathbf{F}(\mathbf{v}).
\end{align}
\end{small}
Here, if shunt elements are non-negligible, the entries in {\small $\mathbf{w}$} are generally not the same.
\end{itemize}
As defined in \eqref{eqn:F}, {\small $\mathbf{F}()$} is the continuously differentiable function that maps any {\small $\mathbf{v}$} into its corresponding {\small $\mathbf{s}$}. We denote the Jacobian of {\small $\mathbf{F}()$} at {\small $\mathbf{v}$} by {\small $\mathbf{J}_\mathbf{F}(\mathbf{v})$}.

In practice, the nodal voltages and branch currents should satisfy certain bounds on their magnitudes. With the above notations, we write these security constraints in \eqref{eqn:sec1}-\eqref{eqn:sec3}, where {\small $V_j^\mathrm{min}$},{\small $V_j^\mathrm{max}$},{\small $I_{jk}^\mathrm{max}$} are pre-specified positive real constants. Note, the branch current constraints are written as \eqref{eqn:sec3} due to \eqref{eqn:bc}.
\begin{small}
\begin{equation} \label{eqn:sec1}
f^\mathrm{V,low}_j(\mathbf{v})\triangleq|v_j|^2-\left(V_j^\mathrm{min}\right)^2>0,~\forall j\in\mathcal{N}^{PQ}.
\end{equation}
\begin{equation} \label{eqn:sec2}
f^\mathrm{V,up}_j(\mathbf{v})\triangleq-|v_j|^2+\left(V_j^\mathrm{max}\right)^2>0,~\forall j\in\mathcal{N}^{PQ}.
\end{equation}
\begin{equation} \label{eqn:sec3}
f^\mathrm{I,branch}_{jk}(\mathbf{v})\triangleq-|a_{jk}v_0+\mathbf{c}_{jk}^T\mathbf{v}|^2+\left(I_{jk}^\mathrm{max}\right)^2>0,~\forall jk\in\mathcal{E}.
\end{equation}
\end{small}
\subsection{$\mathcal{V}$-Control}
We recall definitions and results from \cite{Vcontrol}; they will be frequently used in this paper.
\begin{defn} \label{defn:I}
For a set {\small $\mathcal{V}$} of complex nodal voltages, we define
\begin{itemize}
\item {\small $\mathcal{V}$} is \textit{secured} if \eqref{eqn:sec1}-\eqref{eqn:sec3} are satisfied {\small $\forall\mathbf{v}\in\mathcal{V}$};
\item {\small $\mathcal{V}$} is \textit{a domain of uniqueness} if {\small $\mathbf{F}(\mathbf{v})=\mathbf{F}(\mathbf{v}')\Rightarrow\mathbf{v}=\mathbf{v}',~\forall\mathbf{v},\mathbf{v}'\in\mathcal{V}$};
\item {\small $\mathcal{V}$} is \textit{non-singular} if {\small $\forall\mathbf{v}\in\mathcal{V}$}, the load-flow Jacobian {\small $\mathbf{J}_\mathbf{F}(\mathbf{v})$} is non-singular.
\end{itemize}
Further, for a set {\small $\mathcal{S}$} of nodal power injections, we define
\begin{itemize}
\item {\small $\mathcal{S}$} is \textit{a domain of} {\small $\mathcal{V}$}-\textit{control} if for any continuous path {\small $\mathbf{v}(t):[0,1]\rightarrow\mathbb{C}^{N}$} such that {\small $\mathbf{v}(0)\in\mathcal{V}$} and {\small $\mathbf{F}(\mathbf{v}(t))\in\mathcal{S},\forall t\in[0,1]$}, we have {\small $\mathbf{v}(t)\in\mathcal{V},\forall t\in[0,1]$}.
\end{itemize}
\end{defn}
In Definition \ref{defn:I}, the concept of {\small $\mathcal{V}$}-control can be interpreted as follows: Keep the continuous trajectory {\small $\mathbf{v}(t)$} in {\small $\mathcal{V}$} by maintaining the continuous trajectory {\small $\mathbf{s}(t)$} in {\small $\mathcal{S}$}. On the basis of {\small $\mathcal{V}$}-control, if {\small $\mathcal{V}$} is further secured, then the electrical state is guaranteed to satisfy the security constraints.

Here we might think that, for any {\small $\mathcal{S}$}, the existence plus uniqueness of the load-flow solution in {\small $\mathcal{V}$} are sufficient for {\small $\mathcal{S}$} to be a domain of {\small $\mathcal{V}$}-control. But, this is not true as discussed in \cite{Vcontrol}. 

Below, we recall a theorem from \cite{Vcontrol}, which gives sufficient conditions for {\small $\mathcal{S}$} to be a domain of {\small $\mathcal{V}$}-control. As can be seen, we need the openness and non-singularity of {\small $\mathcal{V}$} in addition to the existence and uniqueness of the load-flow solution.
\begin{thm}[Lemma 2 of \cite{Vcontrol}] \label{thm:AA}
Let {\small $\mathcal{V}$} be a set of complex nodal voltages and {\small $\mathcal{S}$} be a set of nodal power injections. Assume that
\begin{enumerate}
\item {\small $\mathcal{V}$} is open and non-singular;
\item {\small $\forall\mathbf{s}\in\mathcal{S}$}, there is a unique {\small $\mathbf{v}\in\mathcal{V}$} such that {\small $\mathbf{F}(\mathbf{v})=\mathbf{s}$}.
\end{enumerate}
Then there exists a continuous mapping {\small $\mathbf{G}:\mathcal{S}\rightarrow\mathcal{V}$} such that {\small $\mathbf{F}(\mathbf{G}(\mathbf{s}))=\mathbf{s},\forall\mathbf{s}\in\mathcal{S}$}, and {\small $\mathcal{S}$} is a domain of {\small $\mathcal{V}$}-control.
\end{thm}
\subsection{Theoretical Foundations}
The method in this paper uses Theorem \ref{thm:AA}. However, the uniqueness and existence condition (i.e., the second condition) in Theorem \ref{thm:AA} is difficult to verify in practice. To address this issue, we propose two new theorems that give sufficient conditions for uniqueness (Theorem \ref{thm:BB}) and existence (Theorem \ref{thm:CC}); they form the basis for the method proposed in the rest of the paper. The proofs are in Appendix. 
\begin{thm} \label{thm:BB}
If the set {\small $\mathcal{V}$} of complex nodal voltages is non-singular and convex, then it is a domain of uniqueness.
\end{thm}
\begin{thm} \label{thm:CC}
Let {\small $\mathcal{V}$} be a set of complex nodal voltages, {\small $\mathcal{S}$} be a set of nodal power injections, and {\small $\partial\mathcal{V}$} denote the topological boundary of {\small $\mathcal{V}$}. Assume that
\begin{enumerate}
\item {\small $\mathcal{V}$} is bounded, open and non-singular;
\item {\small $\mathcal{S}$} is connected \footnote{{\small $\mathcal{S}$} is connected if {\small $\mathcal{S}$} itself and the empty set are the only subsets that are both closed and open in {\small $\mathcal{S}$}. For $\mathcal{S}$ to be connected, a
sufficient condition is that $\mathcal{S}$ is \emph{path-connected}, i.e., any two
points in $\mathcal{S}$ can be connected by a continuous path in $\mathcal{S}$.};
\item {\small $\mathbf{F}(\mathcal{V})\bigcap\mathcal{S}$} is not empty;
\item {\small $\mathbf{F}(\partial\mathcal{V})\bigcap\mathcal{S}$} is empty.
\end{enumerate}
Then, for any {\small $\mathbf{s}\in\mathcal{S}$}, there exists a {\small $\mathbf{v}\in\mathcal{V}$} such that {\small $\mathbf{F}(\mathbf{v})=\mathbf{s}$}.
\end{thm}
In essence, Theorem \ref{thm:CC} asserts that every {\small $\mathbf{s}$} in {\small $\mathcal{S}$} has a load-flow solution in {\small $\mathcal{V}$}, provided that (i) at least one {\small $\mathbf{s}^\star$} in {\small $\mathcal{S}$} has a load-flow solution in {\small $\mathcal{V}$}, and that (ii) it is impossible for any {\small $\mathbf{s}$} in {\small $\mathcal{S}$} to have a load-flow solution at the boundary of {\small $\mathcal{V}$}. Intuitively, this is because: If there would be an {\small $\mathbf{s}^{\star\star}$} in {\small $\mathcal{S}$} that has no load-flow solution in {\small $\mathcal{V}$}, then in order to move from {\small $\mathbf{s}^\star$} to {\small $\mathbf{s}^{\star\star}$}, the trajectory in the voltage space must either hit a singular point in {\small $\mathcal{V}$} or exit {\small $\mathcal{V}$} by crossing the boundary {\small $\partial\mathcal{V}$}; but this is made impossible by the 1st and the 4th conditions in Theorem \ref{thm:CC}.
\begin{rem}
In the literature, many results have been given with respect to the load-flow solvability (see e.g., \cite{MIlic1,Chiang1,MIlic2,Chiang2,Bolo,IB,EU,EUThree,Solva2,Solva1,JWS,AN}). Different from these results, the proposed Theorem \ref{thm:CC} uses mainly topological properties.
\end{rem}
\begin{rem}
The proof of Theorem \ref{thm:BB} depends only on {\small $\mathbf{F}()$} being quadratic in rectangular representation; and the proof of Theorem \ref{thm:CC} depends only on {\small $\mathbf{F}()$} being differentiable.
\end{rem}

\section{The Admissibility Problem and a Solution Framework}
\subsection{Problem Formulation} \label{sec:pbf}
First, let us define ``admissibility'' in a formal way.
\begin{defn} \label{defn:II}
Given any electrical state {\small $\mathbf{v}^\mathrm{initial}$} and nodal power injection set {\small $\mathcal{S}^\mathrm{uncertain}$} such that
\begin{enumerate}[label = (I\arabic*)]
\item {\small $\mathbf{v}^\mathrm{initial}$} fulfills \eqref{eqn:sec1}-\eqref{eqn:sec3} and {\small $\mathbf{J}_\mathbf{F}(\mathbf{v}^\mathrm{initial})$} is non-singular;
\item {\small $\mathcal{S}^\mathrm{uncertain}$} is compact (i.e., closed and bounded) and includes {\small $\mathbf{F}(\mathbf{v}^\mathrm{initial})$},
\end{enumerate}
we say {\small $\mathcal{S}^\mathrm{uncertain}$} is \textit{admissible} for {\small $\mathbf{v}^\mathrm{initial}$} if, for any continuous function {\small $\mathbf{v}(t),~t\in[0,1]$} that starts at {\small $\mathbf{v}(0)=\mathbf{v}^\mathrm{initial}$}, we have that {\small $\mathbf{v}(t)$} fulfills \eqref{eqn:sec1}-\eqref{eqn:sec3} and {\small $\mathbf{J}_\mathbf{F}(\mathbf{v}(t))$} is non-singular {\small $\forall t\in[0,1]$} as long as {\small $\mathbf{F}(\mathbf{v}(t))\in\mathcal{S}^\mathrm{uncertain},~\forall t\in[0,1]$}.
\end{defn}
In other words, given any initial electrical state {\small $\mathbf{v}^\mathrm{initial}$} and set {\small $\mathcal{S}^\mathrm{uncertain}$} of nodal power injections that satisfy (I1)(I2), if {\small $\mathcal{S}^\mathrm{uncertain}$} is admissible for {\small $\mathbf{v}^\mathrm{initial}$}, then any continuous trajectory of the electrical state is ensured to remain non-singular and fulfill the security constraints. 

Next, we formulate the admissibility problem as follows.
\\
\noindent
\textbf{Admissibility Problem:} Given {\small $\mathbf{v}^\mathrm{initial}$} and {\small $\mathcal{S}^\mathrm{uncertain}$} that satisfy (I1)(I2), is {\small $\mathcal{S}^\mathrm{uncertain}$} admissible for {\small $\mathbf{v}^\mathrm{initial}$}?

\subsection{Solution Framework} \label{sec:frame}
Observe that, by Definition \ref{defn:I}, {\small $\mathcal{S}^\mathrm{uncertain}$} is admissible for {\small $\mathbf{v}^\mathrm{initial}$} if there exists a set {\small $\mathcal{V}$} such that
\begin{enumerate}[label = (O\arabic*)]
\item {\small $\mathcal{V}$} is secured and non-singular;
\item {\small $\mathcal{S}^\mathrm{uncertain}$} is a domain of {\small $\mathcal{V}$}-control;
\item {\small $\mathbf{v}^\mathrm{initial}\in\mathcal{V}$}. \footnote{This last item is required to ensure that the electrical state remains in $\mathcal{V}$, and does not necessarily follow from (I2).}
\end{enumerate}
By this observation, our framework consists in constructing a set {\small $\mathcal{V}$} such that the hypotheses (O1)-(O3) are satisfied. In step 1, we find a large open set {\small $\tilde{\mathcal{V}}$} that is non-singular and convex, using some sufficient conditions on non-singularity (e.g., \cite{NonSingular}). Then, {\small $\mathcal{V}$} is the intersection of {\small $\tilde{\mathcal{V}}$} and the security constraints. According to Theorem \ref{thm:BB}, the obtained set {\small $\mathcal{V}$} is open, secured, non-singular, and a domain of uniqueness. Thus, (O1) is fulfilled.

In step 2, we first verify (O3) by inspection. Then, we test whether {\small $\forall \mathbf{s}\in\mathcal{S}^\mathrm{uncertain}$}, there is no load-flow solution at the boundary {\small $\partial\mathcal{V}$}. This is done by checking the infeasibility of a number of optimization problems. By Theorem \ref{thm:CC}, this will guarantee that there exists a load-flow solution {\small $\mathbf{v}\in\mathcal{V}$} for any {\small $\mathbf{s}\in\mathcal{S}^\mathrm{uncertain}$} (assuming that {\small $\mathcal{S}^\mathrm{uncertain}$} is connected, which can be easily verified). Further, by Theorem \ref{thm:AA}, this will guarantee that (O2) is satisfied.

Specifically, the framework is described below.

\noindent\rule[0.5ex]{\linewidth}{1.5pt}
\begin{LaTeXdescription}
\item [        Framework] $~$
\begin{enumerate}[label=\textbf{(Step \arabic*)},leftmargin=0.8cm]
\item Construct {\small $\mathcal{V}$} as follows:
\begin{itemize}[leftmargin=0.36cm]
\item Find continuous functions {\small $f_\ell()$}, {\small $\ell\in\{1,...,\tilde{L}\}$} such that {\small $\tilde{\mathcal{V}}\triangleq\lc \mathbf{v}: f_\ell(\mathbf{v})>0,~\ell=1,...,\tilde{L}\rc $} is non-singular and convex (e.g., using the conditions in \cite{NonSingular});
\item Let {\small $f_\ell()>0$}, {\small $\ell\in\{\tilde{L}+1,...,L\}$} be the security constraints \eqref{eqn:sec1}-\eqref{eqn:sec3};
\item Then, let {\small $\mathcal{V}\triangleq\{\mathbf{v}:f_\ell(\mathbf{v})>0,~\ell=1,...,L\}$}.
\end{itemize}
\item Test whether
\begin{itemize}[leftmargin=0.36cm]
\item {\small $\mathbf{v}^\mathrm{initial}\in\mathcal{V}$};
\item {\small $\mathcal{S}^\mathrm{uncertain}$} is connected;
\item The following optimization problems are infeasible for all {\small $\ell$}.
\begin{small}
\begin{align}
{\normalsize\mathbf{[P0}({\small \ell})\mathbf{]}\quad\mathrm{min}} &~\sum_{j=1}^N\left(\mathrm{Re}(v_{j})+\mathrm{Im}(v_{j})\right)\nonumber\\
{\normalsize \mathrm{s.t.:}} & ~f_{\ell'}(\mathbf{v})\geq 0,~\forall \ell'\in\{1,...,L\}\setminus\{\ell\}, \nonumber\\
  &~f_\ell(\mathbf{v})=0,\nonumber\\
  &~\mathbf{F}(\mathbf{v})\in\mathcal{S}^\mathrm{uncertain}.\nonumber
\end{align}
\end{small}
\end{itemize}
If all three tests succeed, then declare that {\small $\mathcal{S}^\mathrm{uncertain}$} is admissible for {\small $\mathbf{v}^\mathrm{initial}$}. Otherwise, we are unsure of the admissibility.
\end{enumerate}
\end{LaTeXdescription}
\noindent\rule[0.5ex]{\linewidth}{1.5pt}

\begin{figure}[t]
\centering
\includegraphics[scale=0.66]{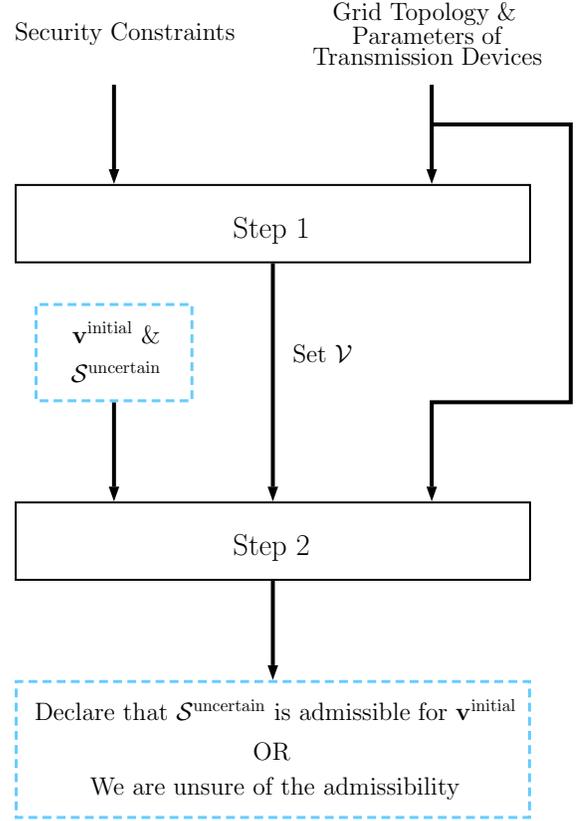}
\caption{Flow chart of the framework.}
\label{fig:FlowChart}
\end{figure}

For this framework, we highlight its structure in Figure \ref{fig:FlowChart} and propose the following theorem on its validity.
\begin{thm} \label{thm:DD}
The above framework is correct in the sense that whenever it declares {\small $\mathcal{S}^\mathrm{uncertain}$} admissible for {\small $\mathbf{v}^\mathrm{initial}$}, it \emph{is} so.
\end{thm}
The proof is in Appendix. In the next section, we develop a polynomial-time method that uses this framework.
\begin{rem} \label{rem:1}
The framework can be extended straightforwardly to multi-phase grids where all non-slack buses operate in {\small $PQ$} mode. Related details can be found in Appendix.
\end{rem}

\section{A Polynomial-Time Method} \label{sec:method}
In this section, we apply the framework to develop a polynomial-time method that, correspondingly, has two steps.
\subsection{Step 1 of the Method}
As in the first step of our framework, we need to concretely find an open set {\small $\tilde{\mathcal{V}}$} that is non-singular and convex. First, take into account that
\begin{itemize}
\item The singularity of load-flow Jacobian usually occurs due to high power generation and consumption;
\item High power generation and consumption are linked to large magnitudes of the branch and nodal currents.
\end{itemize}
Therefore, for {\small $\tilde{\mathcal{V}}$} to be non-singular, we need to ensure that no state in {\small $\tilde{\mathcal{V}}$} has very large branch and nodal currents. Based on this consideration, we let {\small $\tilde{\mathcal{V}}$} be
\begin{small}
\begin{align} \label{eqn:Vtilde}
\tilde{\mathcal{V}}\triangleq\Big\{\mathbf{v}:&\tilde{f}^\mathrm{I,branch}_{jk}(\mathbf{v})>0,~\forall jk\in\mathcal{E}, \nonumber \\
&\tilde{f}^\mathrm{I,node}_{j}(\mathbf{v})>0,~\forall j\in\mathcal{N}^{PQ}
\Big\}
\end{align}
\end{small}
with
\begin{small}
\begin{equation} \label{eqn:sec4}
\tilde{f}^\mathrm{I,branch}_{jk}(\mathbf{v})\triangleq-|a_{jk}v_0+\mathbf{c}_{jk}^T\mathbf{v}|^2+\left(I^\mathrm{branch}_{jk}\right)^2,
\end{equation}
\begin{equation} \label{eqn:sec5}
\tilde{f}^\mathrm{I,node}_{j}(\mathbf{v})\triangleq-|\mathrm{Row}_j(\mathbf{Y}_{LL})(\mathbf{v}-\mathbf{w})|^2+\left(I^\mathrm{node}_{j}\right)^2.
\end{equation}
\end{small}
Here, {\small $\mathrm{Row}_j()$} means the {\small $j$}-th row of a matrix, and {\small $I^\mathrm{branch}_{jk}$}, {\small $I^\mathrm{node}_{j}$} are some auxiliary constants. Obviously, the set {\small $\tilde{\mathcal{V}}$} defined in \eqref{eqn:Vtilde} is already open and convex. Therefore, we need to find appropriate values for constants {\small $I^\mathrm{branch}_{jk}, ~{jk}\in\mathcal{E}$} and {\small $I^\mathrm{node}_{j},~j\in\mathcal{N}^{PQ}$} such that the set {\small $\tilde{\mathcal{V}}$} is non-singular.

To this end, recall that a necessary condition for {\small $\mathbf{J}_\mathbf{F}(\mathbf{v})$} to be singular is given by (8) in \cite{NonSingular} as follows:
\begin{small}
\begin{equation} \label{eqn:necessarycond}
\exists m\in\mathcal{N}^{PQ}{\normalsize \text{ such that }}\sum_{n=1}^{N}|(\mathbf{Y}_{LL}^{-1})_{m,n}(\mathbf{i})_n|\geq|(\mathbf{v})_m|.
\end{equation}
\end{small}
Clearly, if none of the elements in {\small $\tilde{\mathcal{V}}$} satisfies this necessary condition for singularity, then {\small $\tilde{\mathcal{V}}$} is non-singular. Thus, we have the following proposition, where {\small $\|\cdot\|_1$} is the {\small $\ell_1$} norm. The corresponding proof can be found in Appendix.
\begin{prop} \label{prop:1}
The set {\small $\tilde{\mathcal{V}}$} defined in \eqref{eqn:Vtilde} is non-singular if the following optimization problems are infeasible for all {\small $m,n\in\mathcal{N}^{PQ}$} and {\small $\psi,\phi\in\{1,-1\}$}.
\begin{small}
\begin{align}
{\normalsize\mathbf{[P1}({\small m,n,\psi,\phi})\mathbf{]}}&{\normalsize\quad\mathrm{min}}~\sum_{j=1}^N\left(\mathrm{Re}(v_{j})+\mathrm{Im}(v_{j})\right)\nonumber\\
{\normalsize \mathrm{s.t.:~}} &\tilde{f}^\mathrm{I,branch}_{jk}(\mathbf{v})\geq0,~\forall {jk}\in\mathcal{E},\nonumber\\
&\tilde{f}^\mathrm{I,node}_{j}(\mathbf{v})\geq0,~\forall {j}\in\mathcal{N}^{PQ},\nonumber\\
&\|\mathrm{Row}_m(\mathbf{Y}_{LL}^{-1})\|_1\Big(\psi\mathrm{Re}\big(\mathrm{Row}_n(\mathbf{Y}_{LL})(\mathbf{v}-\mathbf{w})\big)\nonumber\\
&+\phi\mathrm{Im}\big(\mathrm{Row}_n(\mathbf{Y}_{LL})(\mathbf{v}-\mathbf{w})\big)\Big)\geq|(\mathbf{v})_m|,\nonumber\\
&\psi\mathrm{Re}\big(\mathrm{Row}_n(\mathbf{Y}_{LL})(\mathbf{v}-\mathbf{w})\big)\geq0,\nonumber\\
&\phi\mathrm{Im}\big(\mathrm{Row}_n(\mathbf{Y}_{LL})(\mathbf{v}-\mathbf{w})\big)\geq0. \nonumber
\end{align}
\end{small}
\end{prop}

By above reasoning, we develop the first step of the method below. 

\noindent\rule[0.5ex]{\linewidth}{1.5pt}
\begin{LaTeXdescription}
\item [        Method (Step 1)] $~$
\begin{enumerate}[label=(1-\alph*)]
\item First, take {\small $I_{jk}^\mathrm{branch}=\beta I_{jk}^\mathrm{max}$}, {\small $\forall jk\in\mathcal{E}$}, where {\small $\beta\in(0,1]$} is a fixed scalar. (Note, large {\small $\beta$} is preferred.)
\item Then, let {\small $I_j^\mathrm{node}=\lambda \hat{I}_j^\mathrm{node}$}, {\small $\forall j\in\mathcal{N}^{PQ}$}, where {\small $\lambda$} is a positive scaling factor and {\small $\hat{I}_j^\mathrm{node}$} is some positive reference value for {\small $I_j^\mathrm{node}$}. (Note, we could let {\small $\hat{I}_j^\mathrm{node}$} be the peak nodal current magnitude at bus {\small $j$} in real-world operation, or simply let all {\small $\hat{I}_j^\mathrm{node}$} be the same.)
\item Next, we start with a small {\small $\lambda$} and gradually increase it by either a fixed ratio or a fixed step size, until (i) P1{\small $(m,n,\psi,\phi)$} is no longer simultaneously infeasible for all {\small $m,n\in\mathcal{N}^{PQ}$} and {\small $\psi,\phi\in\{1,-1\}$}; or (ii) the values {\small $I_j^\mathrm{node}$} are impractically large (e.g., well above {\small $1$} p.u.) .
\item With the penultimate value of {\small $\lambda$}, we obtain {\small $I_j^\mathrm{node}$}, {\small $\forall j\in\mathcal{N}^{PQ}$} and the set {\small $\tilde{\mathcal{V}}$} that is defined in \eqref{eqn:Vtilde}.
\item Last, we let {\small $\mathcal{V}=\{\mathbf{v}\in\tilde{\mathcal{V}}:\mathbf{v}{\normalsize \text{ satisfies }\eqref{eqn:sec1}-\eqref{eqn:sec3}}\}$}.
\end{enumerate}
\end{LaTeXdescription}
\noindent\rule[0.5ex]{\linewidth}{1.5pt}

\subsection{Step 2 of the Method}
According to the second step of the proposed framework, our main task amounts to checking the infeasibility of P0{\small $(\ell)$} for every {\small $\ell\in\{1,...,L\}$} as explained in Section \ref{sec:frame}.

Observe that, for each optimization problem P0{\small $(\ell)$}, we have
\begin{itemize}
\item The objective function is polynomial in {\small $\mathrm{Re}(\mathbf{v})$}, {\small $\mathrm{Im}(\mathbf{v})$};
\item {\small $f_\ell(\mathbf{v}),~\ell\in\{1,...,L\}$} are all polynomial in {\small $\mathrm{Re}(\mathbf{v})$}, {\small $\mathrm{Im}(\mathbf{v})$};
\item {\small $\mathbf{F}(\mathbf{v})$} is a system of polynomials in {\small $\mathrm{Re}(\mathbf{v})$} and {\small $\mathrm{Im}(\mathbf{v})$}.
\end{itemize}
Therefore, the optimization problems P0{\small $(\ell)$} become standard polynomial optimizations if we add the following assumption.
\begin{assu} \label{assu:22}
{\small $\mathcal{S}^\mathrm{uncertain}$} is the Cartesian product of {\small $\mathcal{S}^\mathrm{uncertain}_j,~\forall j\in\mathcal{N}^{PQ}$}, and each {\small $\mathcal{S}^\mathrm{uncertain}_j$} is either a convex polygon or a singleton.
\end{assu}
Note that, under Assumption~\ref{assu:22}, {\small $\mathcal{S}^\mathrm{uncertain}$} is a connected set, as it is path-connected.

Furthermore, note that these polynomial optimization problems are not convex. For this reason, we could apply convex relaxation to them and check whether the relaxed problems are infeasible. Indeed, the infeasibility of the relaxed problem implies the infeasibility of the original problem. As proposed in \cite{Lasserre1}, these non-convex polynomial optimization problems can be effectively approximated by a hierarchy of semi-definite programming relaxations. This hierarchy is arranged by a positive integer called relaxation order. As the relaxation order increases, the relaxed problem becomes closer to the original problem, in terms of the optimal value and feasibility. Despite the theoretical beauty of this hierarchy of relaxations, as the number of variables and the relaxation order increase, it gradually becomes computationally intractable. To cope with this issue, a sparsity-exploiting counterpart of this hierarchy is developed later in \cite{Lasserre2,Lasserre3}, where the level of sparsity depends mainly on the cross terms in the polynomial constraints. In \cite{AppMoment1,AppMoment2,AppMoment3}, very nice examples can be found concerning the application of these hierarchies to power systems.

Taking the above into consideration, we develop the second step of the method below. 

\noindent\rule[0.5ex]{\linewidth}{1.5pt}
\begin{LaTeXdescription}
\item [        Method (Step 2)] $~$
\begin{enumerate}[label=(2-\alph*)]
\item Given the set {\small $\mathcal{V}$} obtained in Method (Step 1), check whether {\small $\mathbf{v}^\mathrm{initial}\in\mathcal{V}$}.
\item With the same {\small $\mathcal{V}$} and the sparsity-exploiting hierarchy of semi-definite programming relaxations in \cite{Lasserre2,Lasserre3}, check whether the relaxed P0{\small $(\ell)$} are all infeasible for some relaxation order. (Note, under Assumption \ref{assu:22}, an empirically good choice of the relaxation order is {\small $2$}.)
\item If both (2-a) and (2-b) are true, then we declare that {\small $\mathcal{S}^\mathrm{uncertain}$} is admissible for {\small $\mathbf{v}^\mathrm{initial}$}. Otherwise, we are unsure of the admissibility.
\end{enumerate}
\end{LaTeXdescription}
\noindent\rule[0.5ex]{\linewidth}{1.5pt}
\begin{rem}
A brief description of the sparsity-exploiting hierarchy of semi-definite programming relaxations can be found in Appendix. Moreover, in Appendix, we explain why we require {\small $\mathcal{S}^\mathrm{uncertain}$} to be the Cartesian product of polygons (in Assumption \ref{assu:22}) rather than other convex sets.
\end{rem}

\subsection{Computational Complexity}
We give below a theorem on the computational complexity of the method. Its proof can be found in Appendix.
\begin{thm} \label{thm:EE}
Under Assumption \ref{assu:22}, the proposed method has a polynomial-time complexity.
\end{thm}

\subsection{Implementation Issues}
\begin{enumerate}
\item For a given grid configuration (i.e., topology, line parameters, etc.), the first step of the proposed method needs to be implemented only once.
\item In the first step of the proposed method, the infeasibility of each P1{\small $(m,n,\psi,\phi)$} can be checked independently. Thus, the first step of the method can be implemented in parallel through a multi-core CPU/GPU or a networked computing infrastructure; this is of significance for relatively large {\small $N$}. Similarly, in the second step of our method, the infeasibility of each relaxed P0{\small $(\ell)$} can also be checked independently. Therefore, the second step of the method can be implemented in parallel as well, which means that the proposed method can be deployed for online applications.
\end{enumerate}

\section{Numerical Evaluations} \label{sec:numev}
\begin{figure}[t]
\centering
\includegraphics[scale=0.6]{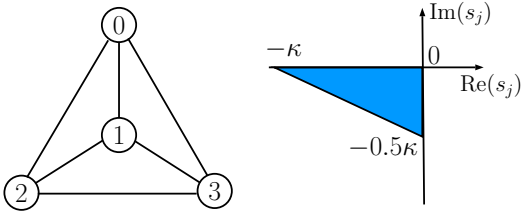}
\caption{Network topology and uncertainty set (in p.u.), where $\kappa$ is a positive real parameter. Note that negative $\mathrm{Re}(s_j),\mathrm{Im}(s_j)$ stand for consumption.}
\label{fig:Meshed}
\end{figure}
\begin{figure}[t]
\begin{center}
\includegraphics[scale=0.6]{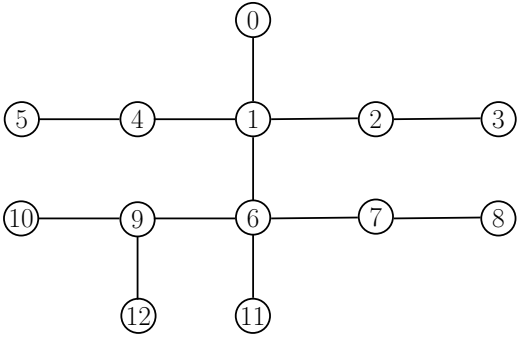}
\caption{Topology of the IEEE 13-Bus Test Feeder.}
\label{fig:IEEE13Bus}
\end{center}
\end{figure}
\begin{figure}[t]
\begin{center}
\includegraphics[scale=0.6]{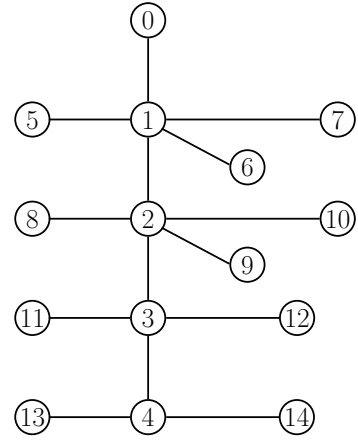}
\caption{Topology of the CIGRE North American LV Distribution Network (residential part).}
\label{fig:CIGRE15Bus}
\end{center}
\end{figure}
In this section, we evaluate the performance of the proposed method in Section \ref{sec:method}, using one meshed grid and two benchmark radial grids in \cite{NumExp,NumExp2,NumExp3}. Topologies of these grids are shown in Figure \ref{fig:Meshed}, \ref{fig:IEEE13Bus} and \ref{fig:CIGRE15Bus}, respectively. For all examples, we assume that (i) the slack-bus voltage is {\small $1$} p.u.; (ii) the relaxation order in the second step of our method is {\small $2$}.

Here, we note that the results in this section are generated on a Macbook Pro, which is equipped with a 2.7 GHz Intel Core i5 CPU and 16 GB 1867 MHz DDR3 memory. In particular, we implement the method using MATLAB tools YALMIP, Mosek and SparsePOP \cite{YALMIP,MOSEK,SparsePOP}.

\subsection{Example 1}
In this example, we consider the meshed grid shown on the left-hand side of Figure \ref{fig:Meshed}. Assume that (i) each transmission line has a series admittance {\small $5-\jmath3.6$} p.u., (ii) the security bounds on nodal voltage magnitudes are {\small $0.95$} and {\small $1.05$} p.u., (iii) the security bounds on branch current magnitudes are {\small $0.6$} p.u., and (iv) {\small $\forall j\in\mathcal{N}^{PQ}$}, {\small $s_j$} belongs to the triangular region on the right-hand side of Figure \ref{fig:Meshed} that specifies {\small $\mathcal{S}^\mathrm{uncertain}$}. Clearly, this grid is stressed when parameter {\small $\kappa\in(0,\infty)$} increases.

Now, let {\small $\mathbf{s}^\mathrm{initial}=\mathbf{0}$} and {\small $\mathbf{v}^\mathrm{initial}=\mathbf{w}$}. We would like to find the maximum value for {\small $\kappa$} such that {\small $\mathcal{S}^\mathrm{uncertain}$} is admissible for {\small $\mathbf{v}^\mathrm{initial}$}. In the first step of the proposed method, we could take {\small $I_{jk}^\mathrm{branch}=I_{jk}^\mathrm{max}$}, {\small $\forall jk\in\mathcal{E}$} for simplicity. Correspondingly, we can choose {\small $I_j^\mathrm{node}=0.8$} p.u., {\small $j\in\mathcal{N}^{PQ}$} and obtain a valid {\small $\mathcal{V}$}. Then, in the second step of the proposed method, we verify that {\small $\mathbf{v}^\mathrm{initial}\in\mathcal{V}$} and find that the maximum value for {\small $\kappa$} to preserve admissibility is {\small $0.35$}. With {\small $\kappa=0.35$}, if {\small $s_j=-\kappa$} p.u., {\small $\forall j\in\mathcal{N}^{PQ}$}, we find that there is a secured and non-singular load-flow solution that has the following features:
\begin{itemize}
\item All nodal voltage magnitudes are {\small $0.9506$} p.u., which indicates that the proposed method is tight in terms of the obtained maximum value of {\small $\kappa$};
\item All branch current magnitudes are much lower than the security bounds;
\item {\small $\forall j\in\mathcal{N}^{PQ}$}, {\small $|i_j|$} is far below {\small $I_j^\mathrm{node}$}, which means that {\small $I_j^\mathrm{node}$} does not limit the performance.
\end{itemize}

\begin{rem}
In this example, we intentionally choose the triangular shape to demonstrate that our method works for polygonal uncertainty sets.
\end{rem}

\subsection{Example 2}
The IEEE 13-Bus Test Feeder is a medium-voltage multi-phase grid, which has shunt elements and a MV/LV transformer \footnote{This transformer lies between buses 2 and 3. By \cite{PowerModel} and \cite{ComAdm}, it is modeled as the serial combination of a winding admittance and an ideal MV/LV transformer. In detail, (i) the winding admittance is positioned between bus 2 and the MV side of the ideal transformer, (ii) bus 3 is directly connected to the LV side of the ideal transformer. In this example, we describe the voltage and current at bus 3 by their equivalents at the MV side of the ideal transformer. Moreover, we do not consider any power limit of this transformer.}. In order to obtain a single-phase grid, we alter the multi-phase IEEE 13-Bus Test Feeder. Specifically,
\begin{itemize}
\item We take the positive-sequence parameters of line configuration 602, \footnote{For each $3\times3$ parameter matrix, we first replace the diagonal positions by their average. Then, we replace all the off-diagonal positions by their average. In this way, the line becomes perfectly transposed, which is characterized by symmetric parameter matrices. As a result, three symmetrical components can be computed without any mutual coupling in the sequence space.} and assume that all the lines are characterized by these parameters;
\item Similarly, we take the positive-sequence equivalent of the transformer;
\item The regulator between buses {\small $0$} and {\small $1$} is removed, as we do not control it.
\end{itemize}
After alteration, the {\small $R/X$} ratio of each transmission line in the resulted single-phase grid is around {\small $0.8$}.

Now, let {\small $\mathbf{s}^\mathrm{initial}$} be the initial nodal power injection, for which {\small $s_j^\mathrm{initial}$} is the average of the IEEE multi-phase benchmark powers at bus {\small $j$}. In addition, let {\small $\mathbf{v}^\mathrm{initial}$} be its high-voltage load-flow solution, which is guaranteed to be unique around {\small $\mathbf{w}$} by theories in \cite{EU}. To ensure that {\small $\mathbf{v}^\mathrm{initial}$} satisfies the security constraints, we choose
\begin{itemize}
\item {\small $V_j^\mathrm{min}=0.9$} p.u. {\small $\forall j\in\mathcal{N}^{PQ}$};
\item {\small $V_j^\mathrm{max}=1.1$} p.u. {\small $\forall j\in\mathcal{N}^{PQ}$};
\item {\small $I_{01}^\mathrm{max}=I_{10}^\mathrm{max}=1$} p.u., {\small $I_{16}^\mathrm{max}=I_{61}^\mathrm{max}=0.45$} p.u., and {\small $I_{jk}^\mathrm{max}=0.3$} p.u. {\small $\forall jk\in\mathcal{E}\setminus\{01,10,16,61\}$}.
\end{itemize}
Assume that, in this medium-voltage grid, the power demands fluctuate significantly. More precisely, {\small $\mathcal{S}^\mathrm{uncertain}$} is a set such that {\small $\forall j\in\mathcal{N}^{PQ},~\mathcal{S}^\mathrm{uncertain}_j=[\kappa\mathrm{Re}(s_j^\mathrm{initial}),0]\times[\kappa\mathrm{Im}(s_j^\mathrm{initial}),0]$}, where {\small $\kappa\in[1,\infty)$} is a scalar. Here, we note that the active and reactive nodal power injections are negative, as there is only power consumption. Evidently, as {\small $\kappa$} increases, {\small $\mathcal{S}^\mathrm{uncertain}$} will eventually fail in the admissibility test for the given {\small $\mathbf{v}^\mathrm{initial}$}. Hence, in the following, we look for the maximum value of {\small $\kappa$} such that {\small $\mathcal{S}^\mathrm{uncertain}$} is admissible for {\small $\mathbf{v}^\mathrm{initial}$}.

According to our method, let us first find proper values for {\small $I^\mathrm{branch}_{jk}, ~{jk}\in\mathcal{E}$} and {\small $I^\mathrm{node}_{j},~j\in\mathcal{N}^{PQ}$}, so that P1{\small $(m,n,\psi,\phi)$} is infeasible simultaneously for all {\small $m,n\in\mathcal{N}^{PQ}$} and {\small $\psi,\phi\in\{1,-1\}$}. For simplicity, we take {\small $I^\mathrm{branch}_{jk}=I_{jk}^\mathrm{max}, ~{jk}\in\mathcal{E}$}. With these {\small $I^\mathrm{branch}_{jk}, ~{jk}\in\mathcal{E}$}, one choice of {\small $I^\mathrm{node}_{j},~j\in\mathcal{N}^{PQ}$} is: {\small $I^\mathrm{node}_1=I^\mathrm{node}_6=0.2$} p.u., {\small $I^\mathrm{node}_2=I^\mathrm{node}_3=I^\mathrm{node}_4=I^\mathrm{node}_5=I^\mathrm{node}_8=0.15$} p.u., and {\small $I^\mathrm{node}_7=I^\mathrm{node}_9=I^\mathrm{node}_{10}=I^\mathrm{node}_{11}=I^\mathrm{node}_{12}=0.1$} p.u.

So far, we have obtained a set {\small $\mathcal{V}$}. Next, using the second step of our method, we find that {\small $\mathbf{v}^\mathrm{initial}\in\mathcal{V}$} and the maximum value for {\small $\kappa$} to preserve admissibility is {\small $1.96$}. When {\small $\kappa=1.96$}, we find that there is a secured and non-singular load-flow solution to {\small $\kappa\mathbf{s}^\mathrm{initial}$}, which has the following features:
\begin{itemize}
\item The lowest nodal voltage magnitude is {\small $|v_8|=0.9016$} p.u.;
\item {\small $|i_{16}|\approx|i_{61}|=0.4128$} p.u., and all the other branch current magnitudes are far below the security bounds;
\item  {\small $|i_j|<I^\mathrm{node}_{j},~\forall j\in\mathcal{N}^{PQ}$}.
\end{itemize}
Thus, in this example, our method is tight in the sense that it almost finds the largest possible value for {\small $\kappa$}. In addition to this tightness, another positive side of our method is the polynomial-time complexity. Specifically,
\begin{itemize}
\item In the first step of the method, the infeasibility of each P1{\small $(m,n,\psi,\phi)$} can be checked in less than {\small $1$} second. And this would be the total execution time if we parallelly check the infeasibility for all P1{\small $(m,n,\psi,\phi)$}. If we sequentially check the infeasibility for all P1{\small $(m,n,\psi,\phi)$}, then the accumulated execution time is {\small $9$} minutes;
\item In the second step of our method, the infeasibility of each relaxed P0{\small $(\ell)$} can be checked in {\small $4-9$} seconds. And this would be the total execution time if we parallelly check the infeasibility for all relaxed P0{\small $(\ell)$}. If we sequentially check the infeasibility for all relaxed P0{\small $(\ell)$}, then the accumulated execution time is around {\small $6$} minutes.
\end{itemize}

\subsection{Example 3}
The residential part of the CIGRE North American LV Distribution Network is a low-voltage split-phase single-phase grid, where every bus is either on the main lateral or directly linked to the main lateral. In this grid, {\small $|i_{jk}|=|i_{kj}|$} holds everywhere, since shunt elements are completely ignored due to short transmission lines. Additionally, the {\small $R/X$} ratios throughout the grid are much larger than {\small $1$}.

We assume that each of the buses {\small $1$-$4$} has an extra energy source. Moreover,
\begin{itemize}
\item Each of these sources is balanced across the neutral line;
\item Each of these sources has an active power generation in {\small $[(1-\kappa)\times20,(1+\kappa)\times20]$} kW, where scalar {\small $\kappa\in[0,1)$};
\item These sources are independent of each other.
\end{itemize}
By fixing the benchmark peak power for the other buses, we construct a set {\small $\mathcal{S}^\mathrm{uncertain}$}. Now, let (i) {\small $\mathbf{s}^\mathrm{initial}$} be the central point in {\small $\mathcal{S}^\mathrm{uncertain}$}, and (ii) {\small $\mathbf{v}^\mathrm{initial}$} be its high-voltage load-flow solution that is guaranteed to be unique around {\small $\mathbf{w}$} by theories in \cite{EU}. To ensure that {\small $\mathbf{v}^\mathrm{initial}$} satisfies the security constraints, we choose
\begin{itemize}
\item {\small $V_j^\mathrm{min}=0.95$} p.u. {\small $\forall j\in\mathcal{N}^{PQ}$};
\item {\small $V_j^\mathrm{max}=1.05$} p.u. {\small $\forall j\in\mathcal{N}^{PQ}$};
\item {\small $I_{01}^\mathrm{max}=1$} p.u., {\small $I_{12}^\mathrm{max}=0.8$} p.u., {\small $I_{23}^\mathrm{max}=0.6$} p.u., {\small $I_{34}^\mathrm{max}=0.5$} p.u., and {\small $I_{jk}^\mathrm{max}=0.4$} p.u. {\small $\forall jk\in\mathcal{E}\setminus\{01,12,23,34\}$}.
\end{itemize}
Similarly to the last example, we look for the maximum value of {\small $\kappa$} such that {\small $\mathcal{S}^\mathrm{uncertain}$} is admissible for {\small $\mathbf{v}^\mathrm{initial}$}. In the first step of our method, we take {\small $I^\mathrm{branch}_{jk}=I_{jk}^\mathrm{max}, ~{jk}\in\mathcal{E}$}. With these {\small $I^\mathrm{branch}_{jk}, ~{jk}\in\mathcal{E}$}, one valid choice of {\small $I^\mathrm{node}_{j},~ j\in\mathcal{N}^{PQ}$} is: {\small $I^\mathrm{node}_j=0.6$} p.u. {\small $\forall j\in\mathcal{N}^{PQ}$}. In this way, we have obtained a set {\small $\mathcal{V}$}. Next, using the second step of our method, we find that {\small $\mathbf{v}^\mathrm{initial}\in\mathcal{V}$} and the maximum value for {\small $\kappa$} to preserve admissibility is {\small $0.11$}. When {\small $\kappa=0.11$}, if each of the four extra sources at buses {\small $1$-$4$} has a power generation {\small $(1+\kappa)\times20$} kW, then there is a secured and non-singular load-flow solution that has the following features:
\begin{itemize}
\item All the nodal voltage magnitudes are close to {\small $1$} p.u.;
\item {\small $|i_{01}|=0.9498$} p.u., {\small $|i_{12}|=0.7813$} p.u., {\small $|i_{23}|=0.5428$} p.u., and all the other branch current magnitudes are far below the security bounds;
\item  {\small $|i_j|<I^\mathrm{node}_{j},~\forall j\in\mathcal{N}^{PQ}$}.
\end{itemize}
Again, our method is tight in the sense that it almost finds the largest possible value for {\small $\kappa$}. In terms of the time cost, we have
\begin{itemize}
\item In the first step of the method, the infeasibility of each P1{\small $(m,n,\psi,\phi)$} can be checked in roughly {\small $1$} second. And this would be the total execution time if we parallelly check the infeasibility for all P1{\small $(m,n,\psi,\phi)$}. If we sequentially check the infeasibility for all P1{\small $(m,n,\psi,\phi)$}, then the accumulated execution time is {\small $12-13$} minutes;
\item In the second step of our method, the infeasibility of each relaxed P0{\small $(\ell)$} can be checked in {\small $7-15$} seconds. And this would be the total execution time if we parallelly check the infeasibility for all relaxed P0{\small $(\ell)$}. If we sequentially check the infeasibility for all relaxed P0{\small $(\ell)$}, then the accumulated execution time is about {\small $10$} minutes.
\end{itemize}

\section{Conclusions}
We have studied the admissibility problem in single-phase microgrids, where the electrical state is represented by complex nodal voltages and controlled by the nodal power injections. In order to test the admissibility, we have developed a framework of solution method, using the recently proposed {\small $\mathcal{V}$}-control. For the theoretical foundation of the framework, we show that if a set {\small $\mathcal{V}$} of complex nodal voltages is non-singular and convex, then it is a domain of uniqueness. In addition, given any set {\small $\mathcal{S}$} of nodal power injections, we have presented topological conditions on {\small $\mathcal{V}$} and {\small $\mathcal{S}$} to guarantee that every element in {\small $\mathcal{S}$} has a load-flow solution in {\small $\mathcal{V}$}, based on the impossibility of obtaining load-flow solutions at the boundary of {\small $\mathcal{V}$}. Within the developed framework, we have established a polynomial-time method that mainly involves the infeasibility check of convex optimizations. The method has been then evaluated on IEEE and CIGRE test grids. Numerical results demonstrate that the method has potential for real-world applications.

\appendices
\section*{Appendix}
\subsection{Proof of Theorem \ref{thm:BB}} \label{apd:BBproof}
\begin{proof}
We prove by contradiction. Let set {\small $\mathcal{V}$} be non-singular and convex. In addition, suppose that there exist {\small $\mathbf{v},\mathbf{v}'\in\mathcal{V}$} such that (i) {\small $\mathbf{F}(\mathbf{v})=\mathbf{F}(\mathbf{v}')$}, and (ii) {\small $\mathbf{v}\neq\mathbf{v}'$}. Owing to convexity, there is {\small $\frac{\mathbf{v}+\mathbf{v}'}{2}\in\mathcal{V}$}. Furthermore, according to the non-singularity in {\small $\mathcal{V}$}, {\small $\mathbf{J}_\mathbf{F}(\frac{\mathbf{v}+\mathbf{v}'}{2})$} is non-singular. However, by the Property 1 in \cite{DomainUnique}, {\small $\mathbf{J}_\mathbf{F}(\frac{\mathbf{v}+\mathbf{v}'}{2})$} should be singular since {\small $\mathbf{F}()$} is a collection of quadratics in rectangular representation. This creates a contradiction and completes the proof.
\end{proof}
\subsection{Proof of Theorem \ref{thm:CC}} \label{apd:CCproof}
\begin{proof}
We need to show that {\small $\mathcal{S}\subseteq\mathbf{F}(\mathcal{V})$}, i.e., {\small $\mathbf{F}(\mathcal{V})\bigcap\mathcal{S}=\mathcal{S}$}. Since {\small $\mathcal{S}$} is connected, its closed and open subsets are {\small $\mathcal{S}$} and the empty set. Based on this, consider that (i) {\small $\mathbf{F}(\mathcal{V})\bigcap\mathcal{S}\subseteq\mathcal{S}$}, and (ii) {\small $\mathbf{F}(\mathcal{V})\bigcap\mathcal{S}$} is not empty, we can prove {\small $\mathbf{F}(\mathcal{V})\bigcap\mathcal{S}=\mathcal{S}$} by showing {\small $\mathbf{F}(\mathcal{V})\bigcap\mathcal{S}$} is both closed and open in {\small $\mathcal{S}$}.

First, the openness of {\small $\mathcal{V}$} implies {\small $\partial \mathcal{V} = \mathrm{cl}(\mathcal{V})\setminus\mathcal{V}$}, where {\small $\mathrm{cl}(\mathcal{V})$} is the closure of {\small $\mathcal{V}$}. Thus, {\small $\mathbf{F}(\mathrm{cl}(\mathcal{V})\setminus\mathcal{V})\bigcap\mathcal{S}$} is empty.
As {\small $\mathcal{V}$} is bounded, we have that the closure {\small $\mathrm{cl}(\mathcal{V})$} is compact. Therefore, by continuity of {\small $\mathbf{F}()$}, {\small $\mathbf{F}(\mathrm{cl}(\mathcal{V}))$} is compact and {\small $\mathbf{F}(\mathrm{cl}(\mathcal{V}))\bigcap\mathcal{S}$} is closed in {\small $\mathcal{S}$}. Since {\small $\mathbf{F}(\mathrm{cl}(\mathcal{V})\setminus\mathcal{V})\bigcap\mathcal{S}$} is empty, we have {\small $\mathbf{F}(\mathrm{cl}(\mathcal{V}))\bigcap\mathcal{S}=\mathbf{F}(\mathcal{V})\bigcap\mathcal{S}$}. So, {\small $\mathbf{F}(\mathcal{V})\bigcap\mathcal{S}$} is closed in {\small $\mathcal{S}$}.

Second, remember that {\small $\mathcal{V}$} is open and non-singular. By the Inverse Function Theorem \cite{InverseFunction}, {\small $\mathbf{F}(\mathcal{V})\bigcap\mathcal{S}$} is open in {\small $\mathcal{S}$}. Thus, {\small $\mathbf{F}(\mathcal{V})\bigcap\mathcal{S}$} is a non-empty, closed and open subset in {\small $\mathcal{S}$}, which means that {\small $\mathbf{F}(\mathcal{V})\bigcap\mathcal{S}=\mathcal{S}$} and completes the proof.
\end{proof}

\subsection{Proof of Theorem \ref{thm:DD}} \label{apd:DDproof}
\begin{proof}
We need to show that {\small $\mathcal{S}^\mathrm{uncertain}$} is a domain of {\small $\mathcal{V}$}-control. By Theorem \ref{thm:AA}, as {\small $\mathcal{V}$} is already open and non-singular, we only need to prove that {\small $\forall\mathbf{s}\in\mathcal{S}^\mathrm{uncertain}$}, there is a unique {\small $\mathbf{v}\in\mathcal{V}$} such that {\small $\mathbf{F}(\mathbf{v})=\mathbf{s}$}.

According to Theorem \ref{thm:BB}, {\small $\mathcal{V}$} is a domain of uniqueness because it is included in a non-singular and convex set. In this way, it suffices to show that, for any {\small $\mathbf{s}\in\mathcal{S}^\mathrm{uncertain}$}, there exists a {\small $\mathbf{v}\in\mathcal{V}$} such that {\small $\mathbf{F}(\mathbf{v})=\mathbf{s}$}. For this purpose, we should check the four conditions in Theorem \ref{thm:CC}.

In Theorem \ref{thm:CC}, the 1st condition is automatically satisfied, since the security constraints imply boundedness. Also, the 2nd condition is satisfied. The 3rd condition follows from {\small $\mathbf{F}(\mathbf{v}^\mathrm{initial})\in\mathcal{S}^\mathrm{uncertain}$}. 

Now, let us focus on the 4th condition. Since P0({\small $\ell$}) is infeasible for all {\small $\ell$}, we have that the set
{\small $$\mathbf{F}\Bigg(\bigcup_{\ell=1}^L \bigg\{ \mathbf{v}: f_\ell(\mathbf{v})=0\text{ and }f_{\ell'}(\mathbf{v})\geq0,~\ell'\in\{1,...,L\}\setminus\{\ell\}\bigg\} \Bigg)$$}
has an empty intersection with {\small $\mathcal{S}^\mathrm{uncertain}$}. Therefore, we can complete the proof by showing that the boundary {\small $\partial\mathcal{V}$} is contained in the set 
{\small $$\bigcup_{\ell=1}^L \bigg\{ \mathbf{v}: f_\ell(\mathbf{v})=0\text{ and }f_{\ell'}(\mathbf{v})\geq0,~\ell'\in\{1,...,L\}\setminus\{\ell\}\bigg\}.$$}
Consider that all $f_\ell()$ are continuous and the topological boundary of {\small $\mathcal{V}$} is the set of points in {\small $\mathbb{C}^{N}$} that are both limit points of {\small $\mathcal{V}$} and limit points of the complement of {\small $\mathcal{V}$}. If {\small $\mathbf{v}\in\partial\mathcal{V}$}, then {\small $\mathbf{v}$} is the limit of some infinite sequence {\small $\mathbf{v}^{(n)}\in\mathcal{V}$}, thus {\small $f_\ell(\mathbf{v}^{(n)})>0$} and {\small $f_\ell(\mathbf{v})\geq 0$} for all {\small $\ell\in\{1,...,L\}$}. Also, {\small $\mathbf{v}$} is the limit of some infinite sequence {\small $\mathbf{v}'^{(n)}$} outside {\small $\mathcal{V}$}. Since there are only finitely many inequalities, there must be at least one inequality, say with index {\small $\ell^\star$}, such that {\small $f_{\ell^\star}(\mathbf{v}'^{(n)})\leq 0$} for an infinite number of indexes {\small $n$}. It follows that {\small $f_{\ell^\star}(\mathbf{v})\leq 0$} and thus {\small $f_{\ell^\star}(\mathbf{v})= 0$}.
\end{proof}
\subsection{Extension of the Framework to Multi-Phase Grids} \label{apd:extension3Phase}
Now, consider a multi-phase grid that consists of one slack bus and {\small $N$ $PQ$} buses. For phase {\small $\gamma\in\{a,b,c\}$} at bus {\small $j\in\mathcal{N}$}, we denote the complex phase-to-ground nodal voltage, nodal current and nodal power injection by {\small $v_j^\gamma$}, {\small $i_j^\gamma$} and {\small $s_j^\gamma$}, respectively. In addition, let
\begin{itemize}
\item {\small $\mathbf{v}_j\triangleq(v_j^a,v_j^b,v_j^c)^T$}, {\small $\mathbf{i}_j\triangleq(i_j^a,i_j^b,i_j^c)^T$} and {\small $\mathbf{s}_j\triangleq(s_j^a,s_j^b,s_j^c)^T$};
\item {\small $\mathbf{v}\triangleq(\mathbf{v}_1^T,...,\mathbf{v}_N^T)^T$}, {\small $\mathbf{i}\triangleq(\mathbf{i}_1^T,...,\mathbf{i}_N^T)^T$} and {\small $\mathbf{s}\triangleq(\mathbf{s}_1^T,...,\mathbf{s}_N^T)^T$}.
\end{itemize}
Then, we have that
\begin{itemize}
\item The branch current of phase {\small $\gamma$} from bus {\small $j$} to {\small $k$} can be represented in the same linear form as \eqref{eqn:bc};
\item The relation between {\small $\mathbf{v},\mathbf{i},\mathbf{s}$} can be compactly written in the same way as \eqref{eqn:I}\eqref{eqn:F}, where {\small $\mathbf{Y}_{LL}$} is {\small $3N\times3N$} in size and {\small $\mathbf{w}=-\mathbf{Y}^{-1}_{LL}\mathbf{Y}_{L0}\mathbf{v}_0$}; \footnote{In \cite{EUThree,EUInvY2,EUInvY3}, it is shown that $\mathbf{Y}_{LL}$ is invertible in multi-phase grids.}
\item The security constraints become
\begin{small}
\begin{equation}
f^\mathrm{V,low}_{j,\gamma}(\mathbf{v})\triangleq|v_j^\gamma|^2-\left(V_{j,\gamma}^\mathrm{min}\right)^2>0,
\end{equation}
\begin{equation}
f^\mathrm{V,up}_{j,\gamma}(\mathbf{v})\triangleq-|v_j^\gamma|^2+\left(V_{j,\gamma}^\mathrm{max}\right)^2>0,
\end{equation}
\begin{equation}
f^\mathrm{I,branch}_{jk,\gamma}(\mathbf{v})\triangleq-|\mathbf{a}_{jk,\gamma}^T\mathbf{v}_0+\mathbf{c}_{jk,\gamma}^T\mathbf{v}|^2+\left(I_{jk,\gamma}^\mathrm{max}\right)^2>0,
\end{equation}
\end{small}
for all {\small $\gamma\in\{a,b,c\}$}, {\small $ j\in\mathcal{N}^{PQ}$} and {\small $jk\in\mathcal{E}$};
\item Definitions \ref{defn:I} and \ref{defn:II} automatically extend;
\item Theorems \ref{thm:AA}, \ref{thm:BB} and \ref{thm:CC} apply to multi-phase grids without modification, since they are formulated and proven with no dependence on the number of phases;
\item The objective function in optimization P0({\small $\ell$}) is changed to {\small $\sum_{\gamma\in\{a,b,c\}}\sum_{j=1}^N\left(\mathrm{Re}(v_{j}^\gamma)+\mathrm{Im}(v_{j}^\gamma)\right)$}.
\end{itemize}

\subsection{Proof of Proposition \ref{prop:1}}  \label{apd:Prop1proof}
\begin{proof}
First, let us construct in \eqref{eqn:Vmn} a collection of sets:
\begin{small}
\begin{align} \label{eqn:Vmn}
\mathcal{V}_{m,n}\triangleq&\Big\{\mathbf{v}\in\tilde{\mathcal{V}}:\|\mathrm{Row}_m(\mathbf{Y}_{LL}^{-1})\|_1\Big(\big|\mathrm{Re}\big(\mathrm{Row}_n(\mathbf{Y}_{LL})(\mathbf{v}-\mathbf{w})\big)\big|\nonumber\\
&+\big|\mathrm{Im}\big(\mathrm{Row}_n(\mathbf{Y}_{LL})(\mathbf{v}-\mathbf{w})\big)\big|\Big)\geq|(\mathbf{v})_m|\Big\},
\end{align}
\end{small}
where {\small $m,n\in\mathcal{N}^{PQ}$}. By inspection, we have that {\small $\mathcal{V}_{m,n}$} is empty when P1{\small $(m,n,\psi,\phi)$} is infeasible {\small $\forall \psi,\phi\in\{1,-1\}$}.

Next, we show that when {\small $\mathcal{V}_{m,n}$} is empty for all {\small $m,n\in\mathcal{N}^{PQ}$}, the necessary condition in \eqref{eqn:necessarycond} holds nowhere in {\small $\tilde{\mathcal{V}}$}. Specifically,
\begin{itemize}
\item By triangle inequality, the emptiness of {\small $\mathcal{V}_{m,n}$} implies that the following inequality holds {\small $\forall\mathbf{v}\in\tilde{\mathcal{V}}$}.
\begin{small}
\begin{equation}
\|\mathrm{Row}_m(\mathbf{Y}_{LL}^{-1})\|_1\big|\mathrm{Row}_n(\mathbf{Y}_{LL})(\mathbf{v}-\mathbf{w})\big|<|(\mathbf{v})_m|.
\end{equation}
\end{small}
\item Consequently for each {\small $m\in\mathcal{N}^{PQ}$}, the following inequality holds {\small $\forall\mathbf{v}\in\tilde{\mathcal{V}}$}, where {\small $\|\cdot\|_\infty$} is the {\small $\ell_\infty$} norm.
\begin{small}
\begin{align}
&\|\mathrm{Row}_m(\mathbf{Y}_{LL}^{-1})\|_1\|\mathbf{Y}_{LL}(\mathbf{v}-\mathbf{w})\|_\infty\nonumber\\
=&\|\mathrm{Row}_m(\mathbf{Y}_{LL}^{-1})\|_1\|\mathbf{i}\|_\infty<|(\mathbf{v})_m|.
\end{align}
\end{small}
\item Further, for each {\small $m\in\mathcal{N}^{PQ}$}, the following holds {\small $\forall\mathbf{v}\in\tilde{\mathcal{V}}$}.
\begin{small}
\begin{align}
\sum_{n=1}^{N}|(\mathbf{Y}_{LL}^{-1})_{m,n}(\mathbf{i})_n|\leq\|\mathrm{Row}_m(\mathbf{Y}_{LL}^{-1})\|_1\|\mathbf{i}\|_\infty<|(\mathbf{v})_m|.
\end{align}
\end{small}
\end{itemize}
Thus, the set {\small $\tilde{\mathcal{V}}$} defined in \eqref{eqn:Vtilde} is non-singular if P1{\small $(m,n,\psi,\phi)$} is infeasible for all {\small $m,n\in\mathcal{N}^{PQ}$} and {\small $\psi,\phi\in\{1,-1\}$}.
\end{proof}
\subsection{Proof of Theorem \ref{thm:EE}} \label{apd:EEproof}
\begin{proof}
In the first step of the proposed method, {\small $I_{jk}^\mathrm{branch}$}, {\small $\forall jk\in\mathcal{E}$} are fixed at the beginning. And we only have a limited number of choices for {\small $I_j^\mathrm{node}$}, {\small $j\in\mathcal{N}^{PQ}$}. Therefore, the infeasibility of P1{\small $(m,n,\psi,\phi)$} is checked by a limited number of times. Furthermore, since P1{\small $(m,n,\psi,\phi)$} is convex, its infeasibility can be checked in polynomial time. Thus, the first step of the method has a polynomial-time complexity.

Similarly, in the second step of the proposed method, we have that (i) the infeasibility of the relaxed P0{\small $(\ell)$} is checked by a limited number of times that depends on the grid size, and that (ii) the relaxed P0{\small $(\ell)$} is convex, hence its infeasibility can be checked in polynomial time. Therefore, the second step of the method also has a polynomial-time complexity.
\end{proof}
\subsection{Sparsity-Exploiting Hierarchy of Semi-Definite Programming Relaxations} \label{apd:tutorial}
In the following, we give a brief description of the sparsity-exploiting hierarchy of semi-definite programming relaxations. Our description is based on the tutorial in \cite{Lasserre4}.

Consider following polynomial optimization problem:
\begin{small}
\begin{align}
{\normalsize \mathrm{min}} &~f_0(\mathbf{x})\nonumber\\
{\normalsize \mathrm{s.t.:}} &~f_k(\mathbf{x})\geq0,~k\in\mathcal{K}.\nonumber
\end{align}
\end{small}
where (i) {\small $\mathcal{K}\triangleq\{1,...,K\}$} is an index set, and (ii) {\small $f_0$} and {\small $f_k,~k\in\mathcal{K}$} are all polynomials in {\small $\mathbf{x}\in\mathbb{R}^M$}.

For each polynomial {\small $f$} of {\small $\mathbf{x}$}, we can express it generically as {\small $f(\mathbf{x})=\sum_{\boldsymbol\alpha\in\mathbb{N}^M} c_f(\boldsymbol\alpha)\mathbf{x}^{\boldsymbol\alpha}$} with some {\small $c_f:\mathbb{N}^M\rightarrow\mathbb{R}$}, where {\small $\mathbf{x}^{\boldsymbol\alpha}=x_1^{\alpha_1}\cdots x_M^{\alpha_M}$}. Then, let us define
\begin{itemize}
\item {\small $\mathcal{M}\triangleq\{1,...,M\}$};
\item {\small $\omega_f\triangleq\lceil\mathrm{deg}(f)/2\rceil$}, where {\small $\mathrm{deg}(f)$} is the degree of {\small $f$};
\item {\small ${\mathcal{I}_{f_k}\triangleq\{j\in\mathcal{M}:\exists\boldsymbol\alpha\in\mathbb{N}^M{\normalsize \text{ such that }}\alpha_j>0}{\normalsize \text{ and }}c_{f_k}(\boldsymbol\alpha)\neq0\}$, $\forall k\in\mathcal{K}$} (i.e., {\small $j$} is in {\small $\mathcal{I}_{f_k}$} if {\small $x_j$} explicitly shows up in the polynomial {\small $f_k$});
\item {\small $\mathcal{E}_{f_0}\triangleq\{\{j,\ell\}\subseteq\mathcal{M}:\exists\boldsymbol\alpha\in\mathbb{N}^M{\normalsize \text{ such that }}\alpha_j>0,~\alpha_\ell>0{\normalsize \text{ and }}c_{f_0}(\boldsymbol\alpha)\neq0\}$} (i.e., {\small $\{j,\ell\}$} belongs to {\small $\mathcal{E}_{f_0}$} if {\small $x_j,x_\ell$} explicitly appear together in a monomial of {\small $f_0$});
\item {\small $\mathcal{E}_{f_k}\triangleq\{\{j,\ell\}\subseteq\mathcal{M}:j,\ell\in\mathcal{I}_{f_k}\},\forall k\in\mathcal{K}$};
\item {\small $\mathcal{A}_\omega^\mathcal{C}\triangleq\{\boldsymbol\alpha\in\mathbb{N}^M:\alpha_j=0,\forall j\not\in\mathcal{C}{\normalsize \text{ and }}\sum_{j=1}^M\alpha_j\leq\omega\}$}, where {\small $\mathcal{C}\subseteq\mathcal{M}$};
\item {\small $\boldsymbol\psi(\mathbf{x},\mathcal{A}_\omega^\mathcal{C})$} is a column vector formed by all monomials {\small $\mathbf{x}^{\boldsymbol\alpha}$}, {\small $\boldsymbol\alpha\in\mathcal{A}_\omega^\mathcal{C}$}.
\end{itemize}
To exploit sparsity, we need to first construct a graph with node set {\small $\mathcal{M}$} and edge set {\small $\mathcal{E}_{f_0}\bigcup\cdots\bigcup\mathcal{E}_{f_K}$}. Next, we find a chordal extension of this graph \cite{Chordal}, and denote the maximal cliques of this chordal extension by {\small $\mathcal{C}_r,~r\in\mathcal{R}\triangleq\{1,...,R\}$} with {\small $R$} being the total number of maximal cliques. Clearly, there exists an index mapping {\small $\theta:\mathcal{K}\rightarrow\mathcal{R}$} such that {\small $\mathcal{I}_{f_k}\subseteq\mathcal{C}_{\theta(k)},\forall k\in\mathcal{K}$}.

Now, the original polynomial optimization problem can be equivalently transformed as follows:
\begin{small}
\begin{align}
{\normalsize \mathrm{min}} &~f_0(\mathbf{x}) \nonumber \\ 
{\normalsize \mathrm{s.t.:}} &~\small{\boldsymbol\psi(\mathbf{x},\mathcal{A}_{\omega-\omega_{f_k}}^{\mathcal{C}_{\theta(k)}}) \boldsymbol\psi(\mathbf{x},\mathcal{A}_{\omega-\omega_{f_k}}^{\mathcal{C}_{\theta(k)}})^T f_k(\mathbf{x})\succeq0,~k\in\mathcal{K}}, \nonumber\\
&~\boldsymbol\psi(\mathbf{x},\mathcal{A}_\omega^{\mathcal{C}_r})\boldsymbol\psi(\mathbf{x},\mathcal{A}_\omega^{\mathcal{C}_r})^T\succeq0,~r\in\mathcal{R}.\nonumber
\end{align}
\end{small}
where {\small $\omega\geq\max\{\omega_{f_0},...,\omega_{f_K}\}$} and ``{\small $\succeq0$}'' means positive semi-definite.

Observe that the above formulation can be rewritten as
\begin{small}
\begin{align}
{\normalsize \mathrm{min}} &~\sum_{\boldsymbol\alpha\in\bigcup_{\ell=1}^R\mathcal{A}_{2\omega}^{\mathcal{C}_\ell}}c_{f_0}(\boldsymbol\alpha)\mathbf{x}^{\boldsymbol\alpha} \nonumber \\ 
{\normalsize \mathrm{s.t.:}} &~\sum_{\boldsymbol\alpha\in\bigcup_{\ell=1}^R\mathcal{A}_{2\omega}^{\mathcal{C}_\ell}}\mathbf{L}_k(\boldsymbol\alpha,\omega)\mathbf{x}^{\boldsymbol\alpha}\succeq0,~k\in\mathcal{K}, \nonumber\\
&~\sum_{\boldsymbol\alpha\in\bigcup_{\ell=1}^R\mathcal{A}_{2\omega}^{\mathcal{C}_\ell}}\mathbf{M}_r(\boldsymbol\alpha,\omega)\mathbf{x}^{\boldsymbol\alpha}\succeq0,~r\in\mathcal{R}. \nonumber
\end{align}
\end{small}
for some real symmetric matrices {\small $\mathbf{L}_k(\boldsymbol\alpha,\omega)$} and {\small $\mathbf{M}_r(\boldsymbol\alpha,\omega)$}.

In this way, a semi-definite programming relaxation of the original problem is obtained by replacing each monomial {\small $\mathbf{x}^{\boldsymbol\alpha}$} with a single real variable {\small $y_{\boldsymbol\alpha}$}:
\begin{small}
\begin{align}
{\normalsize \mathrm{min}} &~\sum_{\boldsymbol\alpha\in\bigcup_{\ell=1}^R\mathcal{A}_{2\omega}^{\mathcal{C}_\ell}}c_{f_0}(\boldsymbol\alpha)y_{\boldsymbol\alpha} \nonumber \\ 
{\normalsize \mathrm{s.t.:}} &~\sum_{\boldsymbol\alpha\in\bigcup_{\ell=1}^R\mathcal{A}_{2\omega}^{\mathcal{C}_\ell}}\mathbf{L}_k(\boldsymbol\alpha,\omega)y_{\boldsymbol\alpha}\succeq0,~k\in\mathcal{K}, \nonumber\\
&~\sum_{\boldsymbol\alpha\in\bigcup_{\ell=1}^R\mathcal{A}_{2\omega}^{\mathcal{C}_\ell}}\mathbf{M}_r(\boldsymbol\alpha,\omega)y_{\boldsymbol\alpha}\succeq0,~r\in\mathcal{R}, \nonumber\\
&~y_\mathbf{0}=1.\nonumber
\end{align}
\end{small}
Obviously, by varying {\small $\omega$}, the size of the above semi-definite programming relaxation changes. In the literature, this parameter {\small $\omega$} is referred to as the relaxation order. With {\small $\omega$} being positive integers, we have a hierarchy of semi-definite programming relaxations.

\subsection{Influence of Assumption \ref{assu:22} on Computational Complexity} \label{apd:computationalComplexity}
The detailed complexity of our polynomial-time method is affected by the set {\small $\mathcal{S}^\mathrm{uncertain}$}. In the following, we explain (i) how our Assumption \ref{assu:22} helps control the detailed complexity, and (ii) what will happen if our Assumption \ref{assu:22} is violated. (Note that, we rely on the concepts and notations in Appendix\ref{apd:tutorial}.)

First, let us organize the relaxed optimization problem into the linear matrix inequality form, which is given as follows:
\begin{small}
\begin{align}
{\normalsize \mathrm{min}}&~\sum_{j=1}^d e_j z^\mathrm{LMI}_j \nonumber \\ 
{\normalsize \mathrm{s.t.:}}&~ -\mathbf{D}_0+\sum_{\ell=1}^d \mathbf{D}_\ell z^\mathrm{LMI}_\ell\succeq0.\nonumber
\end{align}
\end{small}
where 
\begin{itemize}
\item {\small $z^\mathrm{LMI}_j,~j\in\{1,...,d\}$} are optimization variables that correspond to the variables {\small $y_{\boldsymbol\alpha}$} in Appendix\ref{apd:tutorial};
\item {\small $e_j,~j\in\{1,...,d\}$} are coefficients that correspond to the coefficients {\small $c_{f_0}(\boldsymbol\alpha)$} in Appendix\ref{apd:tutorial};
\item {\small $\mathbf{D}_\ell,~\ell\in\{0,...,d\}$} are obtained via (i) constants {\small $1,0,-1$}, and (ii) the matrices {\small $\mathbf{L}_k(\boldsymbol\alpha,\omega)$}, {\small $\mathbf{M}_r(\boldsymbol\alpha,\omega)$} in Appendix\ref{apd:tutorial};
\item {\small $d=\mathrm{card}(\bigcup_{\ell=1}^R\mathcal{A}_{2\omega}^{\mathcal{C}_\ell})$} with {\small $\mathrm{card}()$} being the cardinality (here, recall that {\small $\mathcal{C}_\ell$} and {\small $\mathcal{A}_{2\omega}^{\mathcal{C}_\ell}$} are the maximal cliques and sets described in Appendix\ref{apd:tutorial}).
\end{itemize}
In particular, {\small $\mathbf{D}_0,...,\mathbf{D}_d$} are matrices of size {\small $\zeta$}-by-{\small $\zeta$} with {\small $\zeta=\left(2+\sum_{k=1}^K\mathrm{card}(\mathcal{A}_{\omega-\omega_{f_k}}^{\mathcal{C}_{\theta(k)}})+\sum_{r=1}^R\mathrm{card}(\mathcal{A}_{\omega}^{\mathcal{C}_r})\right)$}.

Note that the above formulation has a dual semi-definite programming in the standard equality form shown below:
\begin{small}
\begin{align}
{\normalsize \mathrm{max}}&~<\mathbf{D}_0,\mathbf{Z}^\mathrm{SE}> \nonumber \\ 
{\normalsize \mathrm{s.t.:}}&~<\mathbf{D}_\ell,\mathbf{Z}^\mathrm{SE}>=e_\ell,~\ell\in\{1,...,d\}, \nonumber\\
&~\mathbf{Z}^\mathrm{SE}\succeq0.\nonumber
\end{align}
\end{small}
where {\small $<\cdot,\cdot>$} is the trace inner product and {\small $\mathbf{Z}^\mathrm{SE}$} is the optimization variable.

Next, take into account that the complexity of solving a semi-definite programming by primal-dual interior-point method depends mainly on {\small $d$} and {\small $\zeta$}, we know that the total complexity is lower if
\begin{itemize}
\item The relaxation order {\small $\omega$} is smaller;
\item The maximal cliques are of smaller sizes.
\end{itemize}
By these thoughts, we discuss as follows:
\begin{itemize}
\item If all the constraints specified by {\small $\mathbf{s}\in\mathcal{S}^\mathrm{uncertain}$} are linear in terms of {\small $\mathrm{Re}(\mathbf{s})$} and {\small $\mathrm{Im}(\mathbf{s})$}, then {\small $\mathbf{F}(\mathbf{v})\in\mathcal{S}^\mathrm{uncertain}$} is a collection of constraints that are quadratic in terms of {\small $\mathrm{Re}(\mathbf{v})$} and {\small $\mathrm{Im}(\mathbf{v})$}. In this way, the minimum possible {\small $\omega$} is 1 and the smallest well-performing {\small $\omega$} is 2. 
\item Furthermore, for each bus {\small $j$}, if the constraints on {\small $s_j$} are independent of the other nodal power injections, then these constraints can be expressed in terms of {\small $v_j$} and only those nodal voltages at the neighbouring buses. In this way, whenever {\small $\mathcal{S}^\mathrm{uncertain}$} is a Cartesian product of {\small $\mathcal{S}^\mathrm{uncertain}_j,~\forall j\in\mathcal{N}^{PQ}$}, we have a collection of smallest possible maximal cliques, which are solely decided by the grid topology.
\end{itemize}
Thus, Assumption \ref{assu:22} helps reduce the detailed complexity.

Finally, in the cases where Assumption \ref{assu:22} is violated, the proposed method could still have a polynomial-time complexity. However, it might not be applicable in practice. To see why, we analyze the following example. Suppose that {\small $\mathcal{S}^\mathrm{uncertain}=\{\mathbf{s}:\|\mathbf{s}\|_2\leq S^\mathrm{max}\}$}. Since all the buses are coupled together, there is only one maximal clique, which has the largest possible size. In this case,
\begin{itemize}
\item {\small $d=\binom{M+2\omega}{2\omega}$};
\item {\small $\zeta=2+\sum_{k=1}^K\binom{M+\omega-\omega_{f_k}}{\omega-\omega_{f_k}}+\binom{M+\omega}{\omega}$};
\item The minimum possible {\small $\omega$} is 2 and the smallest well-performing {\small $\omega$} is 3.
\end{itemize}
With a grid of ten {\small $PQ$} buses (i.e., {\small $M=20$}) and {\small $\omega=3$}, the complexity is already unsuitable for practical applications.

\bibliographystyle{IEEEtran}

\bibliography{REF}
\end{document}